\documentclass[11pt]{article}
\usepackage[mathscr]{eucal}

\usepackage{xypic}
\usepackage{amsfonts}
\usepackage{amsmath}
\usepackage{amsthm}
\usepackage{amssymb}
\usepackage{latexsym}

\newcommand{\bom}{{_{\mathbf{\omega}}}}
\newcommand{\prim}{{\mathtt{prim}}}
\newcommand{\ba}{\overline{A}}
\newcommand{\sss}{{\mathscr{S}}}

\newtheorem{theorem}[equation]{Theorem}
\newtheorem{proposition}[equation]{Proposition}
\newtheorem{corollary}[equation]{Corollary}
\newtheorem{definition}[equation]{Definition}

\newtheorem{lemma}[equation]{Lemma}

\newcommand{\dis}{{\displaystyle}}

\newcommand{\beq}{\begin{equation}\label}
\newcommand{\aand}{\quad{\text{\textsl{and}}\quad}}
\newcommand{\eeq}{\end{equation}}

\newcommand{\iso}{{\;\;\stackrel{\sim}{\longrightarrow}\;\;}}

\newcommand{\vi}{${\sf {(i)}}\;$}
\newcommand{\vii}{${\sf {(ii)}}\;$}

\newcommand{\sset}{\subset}

\newcommand{\id}{{{\mathtt {Id}}}}

\newcommand{\into}{\,\,\hookrightarrow\,\,}
\newcommand{\too}{\,\,\longrightarrow\,\,}

\newcommand{\ad}{{\mathtt{{Ad}}}}

\newcommand{\Spec}{{\mathtt{Spec}}}

\newcommand{\eu}{{{\mathsf {eu}}}}
\newcommand{\dr}{{{\mathsf {DR}}}}
\newcommand{\Hom}{{\mathtt{Hom}}}

\newcommand{\trace}{{\mathtt{Trace}}}
\newcommand{\tr}{{{\mathsf {tr}}^{}}}
\newcommand{\htr}{{{\mathsf{\widehat{tr}}}^{}}}
\newcommand{\dd}{{\mathscr{D}}}

\newcommand{\der}{{\mathtt{Der}^{\,}}}
\newcommand{\derwa}{{\mathtt{Der}}_{\!_B}(A,\omega)}
\newcommand{\derpa}{{\mathtt{Der}}_{\!_{\mathscr{P}}}^{^{_\bullet}}(A,\omega)}

\newcommand{\derrep}{{\mathtt{Der}}_{\!_{\omega_{_{\mathsf{Rep}}}}}}

\newcommand{\derp}{\mathtt{Der}_{_{\!{\mathscr{P}}}}}

\newcommand{\Ker}{{\mathtt {Ker}}}
\newcommand{\aut}{{\mathtt{Aut}}}
\newcommand{\tdash}{\mbox{\tiny{-}}}

\newcommand{\balg}{_{_{\!B{\mathsf{{{\tdash}alg}}}}}}
\newcommand{\palg}{_{_{{\mathscr{P}}\mathsf{{\tdash}alg}}}}

\newcommand{\om}{\omega}
\newcommand{\A}{A}
\newcommand{\Om}{\Omega}

\newcommand{\PP}{{\mathscr{P}}}
\newcommand{\I}{{\mathscr{I}}}

\newcommand{\bulp}{^\bullet_{_{\!\PP\!}}}

\def\C{{\mathbb{C}}}

\def\R{{\mathsf{R}}}

\def\Z{{\mathbb{Z}}}

\def\upa{{{\mathcal U}^{\PP\!}\!A}}
\def\pp{_{_{\!\PP}}}

\def\rep{{\mathsf{Rep}}}
\def\Rep{{\mathscr{R}}}
\def\derb{{{\mathtt{Der}}}_{_{\!B}}}

\def\O{{\sf O}}

\def\T{{\mathsf T}}
\def\Tc{{\mathsf{\check{T}}}}
\def\tc{\check{T}}

\def\g{{\mathfrak{g}}}
\def\gln{{\mathfrak{g}\mathfrak{l}}_n(\C)}
\def\gl{{\mathfrak{g}\mathfrak{l}}}

\def\glinfty{{\mathfrak{g}\mathfrak{l}}_\infty}
\def\sg{{\mathfrak{s}\mathfrak{g}}}
\def\sll{{\mathfrak{s}\mathfrak{l}}}

\def\L{{{\mathfrak{L}}}}

\def\k{{\Bbbk}}
\def\s{{\mathbb{S}}}

\def\ccirc{{{}_{^{^\circ}}}}

\newcommand{\vr}{\varrho}

\newcommand{\al}{\alpha}

\newcommand{\End}{{{\mathtt {End}}}}

\def\cm{{\mathbb{M}}}

\def\ccirc{{}_{^\circ}}
\def\cyc{_{_{\sf{cyclic}}}}

\newcommand{\bone}{{{\mathbf{1}}}}
\newcommand{\lef}{{{\tiny{\sf{left}}}}}
\newcommand{\modd}{{{\tiny{\sf{\tdash}mod}}}}

\newcommand{\GL}{{\mathtt{GL}}}

\def\dq{{\overline{Q}}}

\def\rep{{\mathsf{Rep}}}

\def\g{{\mathfrak{g}}}
\def\gln{{\mathfrak{g}\mathfrak{l}}_n(\C)}

\def\gl{{\mathfrak{g}\mathfrak{l}}}

\def\glinfty{{\mathfrak{g}\mathfrak{l}}_\infty}

\def\cm{{\mathbb{M}}}

\begin{document}

\setlength{\parindent}{6mm}
\setlength{\parskip}{3pt plus 5pt minus 0pt}

\centerline{\Large {\bf Non-commutative Symplectic Geometry,}}
\vskip 3mm
\centerline{\Large {\bf Quiver varieties,$\,$ and$\,$ Operads.}}

\vskip 10mm
\centerline{\large {\sc Victor Ginzburg}}
\vskip 2pt

\hfill{\it to Liza}$\hspace{9mm}$\break

\begin{abstract}
Quiver varieties have recently appeared in
 various different areas of Mathematics
such as
 representation theory of Kac-Moody algebras
and quantum groups,
 instantons on 4-manifolds, and 
resolutions Kleinian singularities.
In this paper, we show that many important affine quiver varieties,
e.g., the Calogero-Moser space,
can be imbedded as coadjoint
orbits in the dual of an appropriate infinite dimensional
Lie algebra. In particular, there is an infinitesimally
transitive action of the Lie algebra in question 
 on the quiver variety. Our construction is based on an extension
of Kontsevich's
formalism of `non-commutative Symplectic geometry'. 
We show  that  this formalism acquires
its most adequate and natural 
formulation in the much more general framework  of
$\PP$-geometry, a `non-commutative geometry' for an algebra
over an arbitrary cyclic Koszul operad.
\end{abstract}
\medskip

\centerline{\sf Table of Contents}\smallskip

$\hspace{30mm}$ {\footnotesize \parbox[t]{115mm}{
1.{ $\;$} {\bf Introduction}\newline
2.{ $\;$} {\bf Non-commutative symplectic geometry}\newline
3.{ $\;$} {\bf Lie algebra associated to a quiver}\newline
4.{ $\;$} {\bf Stabilization: infinite dimensional limit}\newline
5.{ $\;$} {\bf The basics of $\mathbf{\PP}$-geometry}\newline
6.{ $\;$} {\bf Symplectic geometry on a free ${\mathbf{\mathscr{P}}}$-algebra.}
}}

\bigskip

\section{Introduction}
\setcounter{equation}{0}

For the reader's convenience we first remind the definition of quiver varieties.
Let $Q$ be a quiver,
that is a finite oriented graph with  vertex set  $I$.
 Let $V=\{V_i\}_{i\in I}$ be a collection
of finite dimensional $\C$-vector spaces. By a 
representation of $Q$ in $V$ we mean an assignment
of a linear map: $V_i\to V_j$, for any pair $i,j\in I$ and 
each oriented edge of $Q$ with tail $i$ and head $j$.
Let $\Rep(Q,V)$ denote the set of all representations of
$Q$ in $V$, which is a $\C$-vector space. 
The group ${\Huge{\sqcap}}_{i\in I}\;
\GL_{_\C}(V_i)$ acts naturally on $\Rep(Q,V)$, by conjugation.
This action clearly factors through $\,G(V):=\left({\Huge{\sqcap}}_i\,
\GL(V_i)\right)/\C^*,\,$
the
quotient by
 the group $\C^*$  imbedded diagonally, as scalar
matrices, into each of groups $\GL(V_i)$.
Let $\,\g(V)=\left(\oplus\gl(V_i)\right)/\C\,$ 
denote the Lie algebra of the group
$G(V)$.

Let $\dq$ be the 
{\it double} of $Q$, the quiver obtained by adding 
a reverse arrow $a^*$, for every (oriented) arrow 
$a\in Q$. For any $V=\{V_i\}_{i\in I}$,
the vector space $\Rep(\dq, V)$ may be identified naturally
with $T^*\Rep(Q,V)= $ {\it cotangent bundle on}
$\Rep(Q,V).\,$ Hence, $\Rep(\dq, V)$ has a canonical symplectic structure.
Furthermore, the $G(V)$-action on $\Rep(\dq, V)$ is Hamiltonian,
and the corresponding moment map
$\mu: \Rep(\dq, V) \too \g(V)^*$
is given by the following formula:
{\small{
\begin{equation}\label{moment_quiver}
\vr\, \mapsto\, \mu(\vr)=
\Big\lbrace{\{\mu(\vr)_i\}_{_{i\in I}}\in \oplus\gl(V_i)
\enspace\Big|\enspace
\mu(\vr)_{_i}=\!\!\sum_{\big\{{{a\in Q}\atop{\text{head}(a)=i}}\big\}}\!\!\!
\vr(a)\cdot \vr(a^*)\;-\;
\!\!\sum_{\big\{{{a\in Q}\atop{\text{tail}(a)=i}}\big\}}\!\!\!
\vr(a^*)\cdot \vr(a)\Big\rbrace}
\end{equation}}}
\noindent
Here and below, we identify $\g(V)^*$ with a
subspace in $\,\oplus\gl(V_i)\,$ by means of the trace
 pairing: $x,y \mapsto \sum_{i\in I}\, \tr(x_i\cdot
y_i)\,.$ Specifically, we have: 
\[\g(V)^*\;\simeq\; \sg(V)\;:=\;
\lbrace\, x=(x_i)_{i\in I}\;\in \oplus\gl(V_i)\quad\big|\quad\sum\nolimits_{i\in I}\, 
\tr(x_i)=0\,\rbrace\,.\]

\noindent
{\bf Example.}\quad
Let $Q$ be the quiver consisting of a single vertex and a single
edge-loop at this vertex. Thus $\dq$ is the quiver with two edge-loops
at the same vertex. Clearly,
giving a representation of $\dq$ in
the vector space $V=\C^n$ amounts to giving
an arbitrary pair of $n\times n$-matrices.
Therefore, we have: $\Rep(\dq,\C^n)= \gl_n\oplus\gl_n$,
and hence: $G(V)={\mathtt{PGL}}_n$.
The moment map (\ref{moment_quiver}) reduces to the map
$\mu: \gl_n\oplus\gl_n\too \g(V)^*=\sll_n,\,$ given
by the formula: $(x,y) \mapsto [x,y].\,$\bigskip

Next,
fix $\O \subset \sg(V)$, a  closed $\ad\, G(V)$-orbit, and assume that the
group $G(V)$ acts freely on the subvariety $\mu^{-1}(\O) \subset
\Rep(\dq, V)$. Then, the orbit space $\Rep_{\O}(\dq, V)\,:=\,
\mu^{-1}(\O)/G(V)$ is an affine
variety, to be called an {\it affine quiver variety}. 
Thus, by definition:
$\,\dis\Rep_\O(\dq, V)$ $:= \Spec\bigl(\C[\Rep(\dq, V)]^{G(V)}/\I^{G(V)}\bigr)\,,$
where
$\I \subset \C[\Rep(\dq, V)]$ stands for the
defining ideal of the subvariety $\mu^{-1}(\O)$,
and we have used that $\C[\Rep(\dq, V)]^{G(V)}/\I^{G(V)}\,=\,
(\C[\Rep(\dq, V)]/\I)^{G(V)}\,,$ due to reductivity of $G(V)$.
If
$\mu^{-1}(\O)$ is smooth then $\Rep_\O(\dq, V)$ is also smooth, and
the symplectic structure on $\Rep(\dq, V)$ induces, via the symplectic
reduction construction, see [GS], a canonical symplectic
structure
on $\Rep_\O(\dq, V)$.

One of the main results of this paper is

\begin{theorem}\label{main_intro} In the above setting, the symplectic
variety $\Rep_\O(\dq, V)$ can be imbedded as a coadjoint orbit in
the dual of $\L(Q),$ an infinite dimensional Lie algebra canonically
attached to the quiver $Q$.
\end{theorem}

It is implicit in the theorem that the  symplectic structure on
$\Rep_\O(\dq, V)$ goes, under the imbedding, into the canonical
 Kirillov-Kostant symplectic structure on the coadjoint orbit.
Note also that the Lie algebra $\L(Q)$ does not depend
on the representation space~$V$.
\medskip

\begin{remark} 
A choice of Hermitian metric 
on $V$ makes $\Rep(\dq, V)$  a flat hyper-K\"ahler space.
An
equivalence: {\it holomorphic symplectic  reduction}
$\;\Leftrightarrow\;$ {\it hyper-K\"ahler reduction},$\;$
see [Hi], gives, for  many orbits $\O$,
a hyper-K\"ahler structure on the quiver variety 
$\Rep_\O(\dq, V)$. 
Recall further that by a well-known result of Kronheimer [Kr],
any coadjoint orbit in a complex reductive Lie algebra
has a hyper-K\"ahler structure. Based on this analogy,
 N. Hitchin  asked if the Calogero-Moser space
(a special case of
quiver variety, see below) is a coadjoint orbit of some
infinite dimensional Lie algebra. Hitchin's question has
been motivated by the recent work of Berest-Wilson [BW],
who constructed
a transitive 
action of $Aut(A_1)$,
the automorphism group of the Weyl algebra,
on the Calogero-Moser
 space. Theorem \ref{main_intro}
gives a positive answer to Hitchin's question
and
 sheds some new light on the Berest-Wilson construction.
\end{remark}\bigskip

\noindent
{\bf {Strategy of the proof of Theorem \ref{main_intro}.}}$\,\,$ 
 The symplectic structure on $\Rep_\O(\dq, V)$ makes the
coordinate ring $\C[\Rep_\O(\dq, V)]$
 an infinite dimensional Lie  algebra with respect to the Poisson bracket.
We will construct a sequence of Lie algebra morphisms:
\begin{equation}\label{3maps}
\L(Q)\stackrel{\psi}{\too} \C[\Rep(\dq, V)]^{G(V)}
\stackrel{{\mathtt{pr}}}{\too} \C[\Rep(\dq, V)]^{G(V)}/\I^{G(V)} =
\C[\Rep_\O(\dq, V)]\;,
\end{equation}
where $\C[\Rep_\O(\dq, V)]$,
 the coordinate ring, is viewed as a Lie algebra 
 with respect to the Poisson bracket arising from the symplectic
structure
on $\Rep(\dq, V)$, and the map ${\mathtt{pr}}$ stands for the
canonical projection. 

Now, for any  affine symplectic
manifold $X$ and any  point $x\in X$, evaluation at $x$ 
gives a linear function on $\C[X]$, whence induces an evaluation
map: $X \stackrel{{\mathtt{ev}}}{\too}
\C[X]^{^\star}$. Note that the vector space $\C[X]^{^\star}$ 
is an (infinite dimensional)
 Poisson manifold
with Kirillov-Kostant bracket. It is immediate from the definitions that
the 
map: $X \to
\C[X]^{^\star}$ is a morphism of Poisson varieties,
i.e., the induced map on 
polynomial functions is a morphism of Poisson algebras.
Since $X$ is smooth and affine, regular functions on $X$ separate
points of $X$ and, moreover, the differentials of
regular functions span tangent spaces at each point of $X$.
 This implies that the  evaluation map is injective, and that
 the infinitesimal Hamiltonian
action of the
Lie algebra $\C[X]$ (with the Poisson bracket) on the  image 
of  the  evaluation map is 
infinitesimally transitive.  Thus, 
the  evaluation imbedding makes  $X$ a 
coadjoint orbit in $\C[X]^{^\star}$.

Applying the considerations above to the symplectic manifold 
$X=\Rep_\O(\dq, V)$, and dualizing the maps in (\ref{3maps}),
one gets a sequence of Poisson morphisms:
$$\Rep_\O(\dq, V) \stackrel{{\mathtt{ev}}}{\into} \C[\Rep_\O(\dq, V)]^\star 
\stackrel{{\mathtt{pr}}^\star}{\too}
\left(\C[\Rep(\dq, V)]^{G(V)}\right)^\star \stackrel{\psi^\star}{\too}
\L(Q)^\star\;.$$
It will be shown later
 that the composite map above is injective, and the
image of $\Rep_\O(\dq, V)$
 is a coadjoint orbit in $\L(Q)^\star$.
Thus, a key step in proving Theorem \ref{main_intro}
is the construction of  Lie algebra map $\psi$ in 
(\ref{3maps}).

We now illustrate our construction of  $\psi$ in a very 
 special
case,  where $Q$ is
 the quiver consisting of a single vertex and
a single edge-loop (see Example above).
To define the Lie algebra $\L(Q)$, it is convenient to
introduce
an auxiliary 2-dimensional symplectic vector space $(E,\om)$ with basis $x,y$
(corresponding to the two loops in $\dq$)
such that $\om(x,y)=1$. For any $p,q\geq 0$,
we define a $\C$-bilinear map $\{\,,\,\}_{\bom}:
E^{\otimes p}\times E^{\otimes q}
\too E^{\otimes (p+q-2)}$
by the formula:
\begin{eqnarray}\label{mybracket}
\{u_1\otimes u_2\otimes \ldots\otimes u_p\;\;,\;\; v_1\otimes v_2\otimes\ldots\otimes 
v_q\}_{\bom}= \hspace{30mm}\\
\sum\limits_{i=1{}_{{}_{}}}^p
 \sum\limits_{j=1{}_{{}_{}}}^q \om(u_i, v_j)
\cdot u_{i+1}\otimes\ldots\otimes u_p\otimes u_1\otimes\ldots\otimes u_{i-1}\otimes
{}^{{}^{}}v_{j+1}\otimes\ldots\otimes v_q\otimes v_1\otimes\ldots\otimes v_{j-1},
\nonumber
\end{eqnarray}
where
$u_1,\ldots,u_p,\,v_1,\ldots,v_q \in E$.
Assembled together, these maps give a bilinear pairing
$\{-,-\}\bom: TE \times TE\too TE,\,$ where 
$TE = \bigoplus_{i\geq 0}\;
E^{\otimes i}\,$ is
the
tensor algebra of $E$.
Let $\,[TE,TE]\,\subset TE\,$ denote
 the $\C$-linear span of the set $\,\{a\cdot b -b\cdot a\}_{a,b \in
 TE}.\,$ 
\begin{proposition}\label{omL_intro}
The pairing $\{,\}\bom$ gives rise to 
a well-defined Lie algebra structure
on the vector space $\,\L(Q):=TE/[TE,TE]$.
\end{proposition}

\noindent
{\bf Remark.}\quad
One of the goals of the paper is to give an interpretation
of the Lie algebra $\,(TE/[TE,TE]\,,\,\{,\}\bom)\,$
as a sort of Poisson algebra associated to an appropriate
`non-commutative' symplectic variety.\bigskip

To complete our  construction we must define a Lie algebra morphism
$\,\psi: \L(Q)=TE/[TE,TE] \to \C[\Rep(\dq, V)]^{G(V)},\,$
see (\ref{3maps}). As we know, for $V=\C^n$ one has:
$\,\Rep(\dq, V)\simeq \gln\oplus\gln,\,$
and
 $G(V)\simeq {\mathtt{PGL}}_n$. It is convenient
to identify the tensor algebra $TE$ with the
 free associative algebra generated by $x,y$.
We define a $\C$-linear map $\,\tr: TE \to  \C[\gl_n\oplus\gl_n]\,$ by
 assigning
to any non-commutative monomial
$f=x^{k_1}\cdot y^{l_1}\cdot x^{k_2} \cdot\ldots \in TE\,$
a polynomial function $\,\tr f\in \C[\gl_n\oplus\gl_n],\,$ given by the formula:
\begin{equation}\label{trace}
\tr f:\, (X,Y) \mapsto \trace(X^{k_1}\cdot Y^{l_1}\cdot X^{k_2} \cdot\ldots)\quad,
\quad X,Y \in \g=\gl_n\,.
\end{equation}
It is clear that $\tr f\in\C[\gl_n\oplus\gl_n]^{\GL_n}$, and 
that $\tr f=0$ if $f\in [TE,TE]$,
by symmetry of the trace. Thus, the assignment: $f\mapsto \tr f$
gives a well-defined linear map
 $\psi: \L(Q)=TE/[TE,TE] \too \C[\gl_n\oplus\gl_n]^{\GL_n}$.
It turns out that this map is a Lie algebra morphism. This
completes our construction, and the outline of the proof of Theorem
\ref{main_intro}\qed
\bigskip

\noindent
{\bf Example: Calogero-Moser space.}\quad
Let $Q$ be the quiver consisting of a single vertex and a single
edge-loop at this vertex, and assume $\dim V=n$, as above. 
Then, $\g(V)={\mathfrak{p}\mathfrak{g}\mathfrak{l}}_n$.
We will be concerned with the coadjoint orbit
$\O \subset \g(V)^*=\sll_n$, formed by all $n\times n$-matrices of the form:
$s\,\mbox{-}\,\id\,$, where $s$ is a
 rank 1  semisimple 
matrix such that $\mbox{\sl Trace}(s) = \mbox{\sl Trace}(\id)=n$.
Thus, $\O$ is a closed $G$-conjugacy class in $\sll_n$, and
it has been
shown in [W] that 
$$
\mu^{-1}(\O)\;=\;\lbrace (X,Y)\in \sll_n\times\sll_n\quad
\big|\quad [X,Y]+\id\;\enspace \mbox{\sl is a
rank one  semisimple matrix}\rbrace\;,
$$
 is a smooth connected algebraic variety,
and the $\ad G$-diagonal action on $\mu^{-1}(\O)$ is free. The reduced
space $\cm:= \mu^{-1}(\O)/G$ is, according to [KKS] (see also [W]),
 nothing but the phase space of the (rational)
Calogero-Moser integrable system. This
is a smooth affine algebraic symplectic manifold.
Thus, Theorem \ref{main_intro} makes $\cm$ a coadjoint orbit in
$(A/[A,A])^*$, where $A=TE=\C\langle x,y\rangle$. This very special
case was the starting point of our analysis.\qed
\bigskip

An earlier version of this
 paper has been greatly motivated by [BW], whose question led me to the
development of  non-commutative geometry in the special
case of the Calogero-Moser space. The results
presented in \S3 below form a natural generalization
of the Calogero-Moser case.
 This
generalization has been found simultaneously and independently
by L. Le Bruyn [LB1] and the author.\medskip

\noindent
{\footnotesize{{\sl Acknowledgements.} {\it I am grateful to Yu. Berest
and
G. Wilson for explaining to me the results of [BW] prior to their
publication. I have benefited from interesting correspondence
with L. Le Bruyn, especially from his letter [LB1]. 
}}}

\section{Non-commutative Symplectic geometry}
\setcounter{equation}{0}
Throughout this paper we will be working
over a ground field $\k$ of characteristic zero,
and write $\otimes=\otimes_\k$. 
We fix a  commutative unital $\k$-algebra $B$, and for any $B$-bimodule
$M$, write $T^j_BM = M\otimes_{_{\!_B}}\ldots
\otimes_{_{\!_B}}M$ ($j$ factors $M$), which is a $B$-bimodule again.

Let
$A$ be
 a unital associative $\k$-algebra  containing the commutative
algebra $B$ as a subalgebra.
Recall that the free differential envelope of 
$A$ over $B$ is a graded vector space
$\Om_{_{\!B}}^\bullet\!A = \bigoplus_{j\geq 0}\;\Om_{_{\!B}}^j\!A,\,$
where $\,\Om_{_{\!B}}^j\!A= 
A\,\bigotimes_B\, T^j_B(A/B)\,$ is the $B$-bimodule
formed by linear combinations of  expressions $a_0\cdot da_1\ldots  da_j
\in A\otimes T^j_B(A/B)$.
Moreover, it is known, cf. [L],  that there is a
$B$-bimodule isomorphism:
$\Om_{_{\!B}}^\bullet\!A \simeq 
\bigoplus_{{j\geq 0}_{\,}}\;T^j_B(\Om_{_{\!B}}^1\!A)\,,$ 
and there is a $B$-bimodule
super-differential $d:
\Om_{_{\!B}}^\bullet\!A\to\Om_{_{\!B}}^{\bullet +1}A,\,$
making 
${}^{^{}}\Om_{_{\!B}}^\bullet\!A$
an associative differential graded algebra.

Given $\alpha\in \Om_{_{\!B}}^i\!A\,,\, \beta \in \Om_{_{\!B}}^j\!A,$ we put:
$[\alpha,\beta] = \alpha\cdot \beta - (-1)^{ij}\beta\cdot\alpha\,,$
and write $[\Om_{_{\!B}}^\bullet\!A\,,\,\Om_{_{\!B}}^\bullet\!A]$ 
for the $B$-linear span
of all such super-commutators.
Following Karoubi [Ka], see also [L, \S2.6],
define the relative
non-commutative de Rham complex of the pair $(A,B)$ as the differential graded vector
space:
$$\dr_{\!_{B}}^\bullet\!A =  
\Om_{\!_{B}}^\bullet\!A/[\Om_{\!_{B}}^\bullet\!A\,,\,\Om_{\!_{B}}^\bullet\!A]
\quad,\quad \dr_{\!_{B}}^\bullet\!A =\bigoplus\nolimits_{j\geq 0}\,\dr_{\!_{B}}^jA\,,
$$
where the differential and the grading are induced from those on $\Om_{\!_{B}}^\bullet\!A$.
Abusing the notation we will write: $a_0\cdot da_1\ldots  da_j
\in \dr_{\!_{B}}^j A$, meaning the corresponding class modulo commutators.
We have: $\dr_{\!_{B}}^0A=A/[A,A],$ and 
$\,H^0(\dr_{\!_{B}}^\bullet\!A)= \ker(\dr_{\!_{B}}^0A \to \dr_{\!_{B}}^1A) = B\,.$

Let $\der_{\!_B}A$ denote the Lie algebra of all $B$-linear derivations
of
$A$.
Given $\theta \in \der_{\!_B}A$ one introduces, following [K2],
 a Lie operator $L_\theta:
\Om_{\!_{B}}^\bullet\!A\to\Om_{\!_{B}}^\bullet\!A$, resp. a contraction operator
$i_\theta:
\Om_{\!_{B}}^\bullet\!A\to\Om_{\!_{B}}^{\bullet -1}A$, as  a derivation,
resp. a super-derivation, of the associative algebra
$\Om_{\!_{B}}^\bullet\!A$ defined on generators by the formulas:
$$
L_\theta(a_0)=\theta(a_0)\enspace,\enspace
L_\theta(da)= d(\theta(a))\quad\mbox{and}\quad
i_\theta(a_0)=0\enspace,\enspace i_\theta(da)= \theta(a)\quad,\enspace
\forall a_0,a\in A.
$$
It is straightforward to verify that 
 the  induced operators on $\dr_{\!_{B}}^\bullet\!A$,
satisfy the following standard commutation relations:
\begin{equation} \label{identities}
L_\theta = i_\theta\ccirc d + d\ccirc  i_\theta\quad,\quad
[L_\theta, i_\gamma] = i_{[\theta, \gamma]}
\quad,\quad
[L_\theta, L_\gamma] = L_{[\theta, \gamma]}
\quad,\quad
i_\theta\ccirc  i_\gamma = -i_\gamma\ccirc  i_\theta,
\end{equation} 
where all the commutation relations but the last one hold already
in $\Om_{\!_{B}}^\bullet\!A$.

Fix $\om \in \dr_{\!_{B}}^2A$, and set $\derwa=\{\theta\in \der
A\;\,|\,\;L_\theta\om=0\}\,.$
Clearly, $\derwa$ is a Lie subalgebra in $\der_{\!_B}A$.
The assignment:
$\theta \mapsto i_\theta\om$ gives a linear map
$i: \der_{\!_B}A\,\to\,\dr_{\!_{B}}^1A$. The 2-form
$\om\in \dr_{\!_{B}}^2A$ is called {\it non-degenerate} provided the map
$i$ is bijective. 

\begin{lemma}\label{Ko}
Let $\om \in \dr_{\!_{B}}^2A$ be a non-degenerate  2-form such that $d\om=0$ in
$\dr_{\!_{B}}^3A$. Then the map: $\theta \mapsto i_\theta\om$
induces a bijection $i: \derwa \iso (\dr_{\!_{B}}^1A)_{\!_{\sf closed}}\,,$
that is:
$\;\;\theta\in \derwa 
\enspace\Longleftrightarrow\enspace d(i_\theta\om)=0$ in $\dr_{\!_{B}}^2A$.
\end{lemma}

{\sl Proof.} Since, $d\om=0$, we have:
$\,L_\theta\om = i_\theta d\om + di_\theta\om
= di_\theta\om\,.$ Hence, 
$\,\theta \in \derwa\;\Longleftrightarrow\; 0=L_\theta\om=d(i_\theta\om).$
 \qed\medskip

By  Lemma \ref{Ko}, one may invert the isomorphism $i$ 
to obtain a linear bijection
$i^{-1}:\; (\dr_{\!_{B}}^1A)_{\!_{\sf closed}}\,\iso\,$
$\derwa$.
Let: $f\mapsto \theta_f$ denote the map given by
 the  composition:
\begin{equation} \label{Psi} 
 A/[A,A] =\dr_{\!_{B}}^0A\;\stackrel{d}{\too}\;
(\dr_{\!_{B}}^1A)_{\!_{\sf exact}}\;\into\;
(\dr_{\!_{B}}^1A)_{\!_{\sf closed}}
\;\stackrel{i^{-1}}{\too}\;\derwa\,.
\end{equation} 
Using the map: $f\mapsto \theta_f$, we
define a Poisson bracket on $A/[A,A]$ by any of the following
equivalent expressions:
\begin{equation} \label{bom}
 \{f, g\}_{\bom} \;:=\;\;\; i_{_{\theta_f}}(i_{_{\theta_g}}\om)\;=\;
i_{_{\theta_f}}(dg)\;=\;-i_{_{\theta_g}}(df)\;=\;
L_{_{\theta_f}}g\;=\; -L_{_{\theta_g}}f\,.
\end{equation}
Here, in the first expression for $\{f, g\}_{\bom}$ we have used the
composite
map: $\,i_{_{\theta_f}}\ccirc i_{_{\theta_g}}:
\dr_{\!_{B}}^2A\to
\dr_{\!_{B}}^1A\to\dr_{\!_{B}}^0A.$ Other equalities, e.g.:
$\,i_{_{\theta_f}}(i_{_{\theta_g}}\om) = L_{_{\theta_f}}g$, 
follow
from the equation $i_{_{\theta_g}}\om =dg$ (which is the definition of
$\theta_g$),  the obvious identity:
$i_{_{\theta_f}}(dg)= L_{_{\theta_f}}g$,
and the last equation in (\ref{identities}).

\begin{theorem}\label{Lie*} The bracket (\ref{bom}) makes
$A/[A,A]$ into a Lie algebra.
\end{theorem}

{\sl Proof.} Skew symmetry of the bracket is immediate from  (\ref{bom}).
We begin the proof of the Jacobi identity by 
observing that, for any $f\in A/[A,A]$ and $\eta\in \derb{A}$
one has:
\begin{equation} \label{pois1}
i_\eta i_{\theta_f}\om=-L_\eta f\,.
\end{equation}
Further, for any $\xi\,,\,\eta_{_1}\,,\,\eta_{_2}\in
\derb{A}$, commutation relations (\ref{identities})  yield:
$$L_\xi i_{\eta_{_1}} i_{\eta_{_2}}\om -
i_{\eta_{_1}} i_{\eta_{_2}}L_\xi\om =
i_{[\xi, \eta_{_1}]}i_{\eta_{_2}}\om
+i_{\eta_{_1}} i_{[\xi, \eta_{_2}]}\om\,.
$$
In the special case: $\xi=\theta_f$, we have: $L_\xi\om=0$,
hence, the above equation reads:
\begin{equation} \label{pois2}
L_{\theta_f}i_{\eta_{_1}} i_{\eta_{_2}}\om=
i_{[\theta_f, \eta_{_1}]}i_{\eta_{_2}}\om
+i_{\eta_{_1}} i_{[\theta_f, \eta_{_2}]}\om\,.
\end{equation}
We now choose $g\in A/[A,A]$
and put $\eta_{_1}=\theta_g$. Formula (\ref{pois1}) shows that the
LHS of (\ref{pois2}) equals: $-L_{\theta_f}i_{\eta_{_2}}i_{\theta_g}\om=
L_{\theta_f}L_{\eta_{_2}}g\,.$ Similarly,
the second summand on the RHS of (\ref{pois2}) equals: 
$$-
i_{[\theta_f, \eta_{_2}]}i_{\theta_g}\om =
-L_{[\theta_f, \eta_{_2}]}g = L_{\eta_{_2}}L_{\theta_f}g
-L_{\theta_f}L_{\eta_{_2}}g \,.
$$
Thus, writing $\eta=\eta_{_2}$, from (\ref{pois2}) we deduce:
$$
i_{[\theta_f, \theta_g]}(i_{\eta}\om)=
L_\eta L_{\theta_f}g= L_\eta \{f,g\}
=i_{\theta_{\{f,g\}}} i_\eta\om
\,,$$
where the last equality is due to
formula (\ref{pois1}). We conclude 
that the derivation $\,\delta:= [\theta_f, \theta_g] - \theta_{\{f,g\}}
\in \derwa$ has the property that,
for any $\eta\in \derb{A}$, one has: $i_\delta i_\eta\om =0.$
The bijection: $\eta\, \longleftrightarrow\,\alpha=
i_\eta\om$
of Lemma \ref{Ko} now implies that, for any 1-form $\alpha$, one has:
$i_\delta\alpha=0$.

To complete the proof, for any $f\,,\,g\,,\,h\in \dr_{\!_{B}}^0A\,,$ we
write
the identity:
$$(L_{\theta_f}\ccirc L_{\theta_g}-L_{\theta_g}\ccirc L_{\theta_f})h=
L_{[\theta_f, \theta_g]}h= (L_{\theta_{\{f,g\}}+\delta})h
=L_{\theta_{\{f,g\}}}h
 +L_\delta h \,.$$
The leftmost commutator here equals:
$\,\{f,\{g,h\}\} - \{g,\{f,h\}\},\,$ 
and the term $L_{\theta_{\{f,g\}}}h$ on the right equals
$\{\{f,g\}\,,\,h\}$, by definition, see (\ref{bom}).
Finally, we have: $\,L_\delta h = i_\delta(dh)=0\,,$
because of the property of $\delta$ established earlier.
Thus, the identity above yields:
$\,\{f,\{g,h\}\} - \{g,\{f,h\}\} =\{\{f,g\},h\}\,,$
and the Theorem is proved.
\qed\bigskip

Assume that
$\,B=\k\oplus\k\oplus\ldots\oplus\k\,$ (direct sum of $p$ copies of the
 ground field). 
For each $i\in \{1,\ldots,p\},$
let $\bone_i\in B$ denote the idempotent corresponding to the $i$-th
direct
summand $\k$. 

Further, let $V$ be a finite dimensional left $B$-module.
Clearly, giving such a $V$ amounts to
giving a collection of finite dimensional $\k$-vector space 
$\{V_i\}_{\,1\leq i\leq p}$, one for each $i$,
such that $V=\bigoplus_i\;V_i,$ and such that
$\bone_i\in B$ acts as the projector
onto the $i$-th
direct  summand. We consider the algebra $\,\End V:=\End_\k V\,$
of $\k$-linear endomorphisms of $V$. The action of $B$ on $V$
makes $\,V^*:=\Hom_\k(V,\k)\,$ a  right  $B$-module, and
gives an algebra imbedding: $B\into \End V$.
 Hence, left and right multiplication
by $B$ make $\End V$ a $B$-bimodule which is canonically
isomorphic to the  $B$-bimodule $V\otimes_{_\k}V^*$. Further,
 the assignment: 
$$f \mapsto 
\bigl(\tr(\bone_1\!\cdot\! f\!\cdot\!\bone_1)\,\,,\,\,
\tr(\bone_2\!\cdot\! f\!\cdot\!\bone_2)\,\,,\ldots,\,\,
\tr(\bone_p\!\cdot\! f\!\cdot\!\bone_p)\bigr)\,\in\, 
\k\oplus\k\oplus\ldots\oplus\k=B\,$$
gives a canonical
$B$-bimodule trace map $\tr: \End V \too B$.\medskip

\noindent
{\bf {Representation functor.}}\quad\quad\quad
Given a finitely generated associative $B$-algebra $A$, 
 let 
$\Hom\balg(A, \End V)$ denote the affine algebraic variety 
of
all associative algebra homomorphisms $\rho: A\to \End V$,
such that $\rho\big|_B=\id_B$.
Let $\rep(A,V):=\k[\Hom\balg(A, \End V)]$ denote the coordinate ring of
$\Hom\balg(A, \End V)$.
The natural action on $\End V$ of the group
${G(V)}=GL_{_B}(V)$ (of $B$-linear automorphisms of $V$)
 by conjugation induces a
${G(V)}$-action on $\Hom\balg(A, \End V)$.
This gives a $G(V)$-action on $\rep(A,V)$ by
algebra automorphisms. 

The tautological
evaluation map: $\,A\times \Hom\balg(A, \End V) \too \End V\,$
assigns to any element $a\in A$ an $\End V$-valued function $\hat{a}$ on
$\Hom\balg(A, \End V)$. Equivalently, this function may be viewed as
an element $\hat{a}\in \bigl(\rep(A,V)\bigotimes_{\!_B} \End V\bigr)^{G(V)}$.
Taking the trace on the second tensor
factor, one obtains a $G(V)$-invariant $\k$-valued function
$\;\tr(\hat{a})\in \bigl(\rep(A,V)\bigotimes_{\!_B} B\bigr)^{G(V)}=
{\rep(A,V)}^{G(V)}\,.$ The assignment: $\,a\mapsto \tr(\hat{a})\,$
clearly vanishes on $[A,A]$ due to the cyclic symmetry of
the trace map. Thus, it descends to a well-defined $B$-linear
map
\begin{equation}\label{rep_map}
\htr: \;\dr_{\!_{B}}^0{\!A}=A/[A,A]\too {\rep(A,V)}^{G(V)}
\quad,\quad a\mapsto \tr(\hat{a})\;.
\end{equation}

\begin{remark} More generally, for any $p\geq 0$,
the assignment:
$\,a_0\cdot da_1\ldots da_p\,\mapsto\,
\tr(\hat{a}_0\cdot d\hat{a}_1\ldots d\hat{a}_p)\,$
gives a well-defined map from $\dr_{\!_{B}}^p{\!A}$
to the space of $G(V)$-invariant regular  $p$-forms 
(in the ordinary sense) on
the algebraic variety $\Hom\balg(A, \End V)$.
\end{remark}

\section{Lie algebra associated to a Quiver}
\setcounter{equation}{0}

Fix $B$, a commutative $\k$-algebra
and $E$, a finite rank projective $B$-bimodule,
i.e. a projective $B\otimes B^{op}$-module. 
The 
space $\,E^\vee:=\Hom_{_{\lef{B}\modd}}(E,B)\,$ has a 
canonical $B$-bimodule structure given by:   
$\,(b_1\varphi b_2)(e)= \varphi(e\cdot b_1)\cdot b_2\,,$
where $b_1,b_2\in B\,,\, e\in E,$ and $\varphi\in E^\vee$.

A $B$-bimodule map $\om: E\otimes_{\!_{_B}}E \to B$ will be referred to
as a $B$-bilinear form on $E$. For such an $\om$,
the assignment: $e\mapsto \om(-\otimes e)$ gives
a $B$-bimodule map $E\to E^\vee$.
We call $\om$ 
{\it non-degenerate} if the latter map is an isomorphism.
If, furthermore, $\om$ is skew-symmetric,
i.e. $\om(x,y)+\om(y,x)=0,$ for any $x,y\in E$, we
will say that $\om$ is a symplectic $B$-form on $E$.
For example, for any finite rank projective $B$-bimodule $V$, 
the bimodule $E=V\bigoplus V^\vee$ carries a canonical
symplectic $B$-form.

Fix a finite dimensional
$B$-bimodule $E$, and let $\,\A= T_{\!_B}E := \bigoplus_{i\geq 0}\;
T^i_{\!_B}E\,$ be the tensor algebra, a graded associative algebra such that
$T^0_{\!_B}E=B$. For each $i >0$, let
$(T^i_{\!_B}E)\cyc$ denote the quotient
of $T^i_{\!_B}E$ by the $B$-sub-bimodule generated by the
elements:
$$x_1 \otimes_{\!_{_B}} x_2 \otimes_{\!_{_B}} \ldots \otimes_{\!_{_B}}
x_i\;
 -\;
x_i\otimes_{\!_{_B}} x_1 \otimes_{\!_{_B}} \ldots \otimes_{\!_{_B}} x_{i-1}\quad,\quad
\forall x_1,\ldots,x_i\in E\,.$$

The following result was obtained independently by L. Le Bruyn [LB1]
and the author.

\begin{lemma}\label{DR(quiver)} \vi The de Rham complex of
$\A=T_{\!_B}E$
 is acyclic, i.e.,
$H^k(\dr_{\!_B}^\bullet\A) = 0$, for  all $k\geq 1$. Furthermore,
 $H^0(\dr_{\!_B}^\bullet\A) =
B$.

\vii We have: $\;\dr_{\!_B}^0(T_{\!_B}E) = (T_{\!_B}E)\cyc\;,\;\;$ and 
$\;\;\;\dr_{\!_B}^1(T_{\!_B}E) = (T_{\!_B}E)\bigotimes_B\;E\;.$
\end{lemma}

{\sl Proof.} To prove (i), 
we  imitate, following Kontsevich [K2],  the classical
proof of the Poincar\'e lemma. To this end, introduce a ($B$-linear)
Euler derivation
$\eu:  T_{\!_B}E \to T_{\!_B}E$ by letting it act on generators
$x\in E=T_{\!_B}^1E\,$ by: $\eu(x)=x$. The induced
map $L_{\eu}: \dr_{\!_B}^\bullet\A \too \dr_{\!_B}^\bullet\A$
is diagonalizable and has non-negative integral eigenvalues.
Cartan's homotopy
formula: $L_{\eu} = d\ccirc i_\eu +i_\eu\ccirc d$ shows that the de Rham
complex
is quasi-isomorphic to the zero eigen-space of the operator
$L_{\eu}$, which is the subspace $B$ sitting in degree 0.
Part (i) follows. Part (ii) is straightforward.\qed\medskip

From now until the end of the section assume that $B=\k^I,$ where $I$ is
a finite set,
and put
 $A:=T_{\!_B}E$, where $(E,\om)$ is a symplectic $B$-bimodule. 
Using the isomorphism: $E\iso E^\vee =\Hom_{_{\lef{B}\modd}}(E,B),\,$
provided by $\om$, one transports the symplectic structure
from $E$ to $E^\vee$. Let $\,\om^\vee=\sum_r \phi_r\otimes\psi_r\in E\otimes E\,$
be the resulting symplectic $B$-form on  $E^\vee$.
It is straightforward to see that $\,
\sum_r d\phi_r\otimes d\psi_r\in \Omega_{_B}^2A\,$
gives a well-defined {\it closed} and {\it non-degenerate} class in $\dr^2_{\!_B}A\,$,
 to be dented $\om_{_\dr}$.
Thus, the general construction (\ref{bom})
yields a Lie bracket $\{\,,\,\}_{_{\om_{_\dr}}}$ on $A/[A,A]$.
\bigskip

\noindent
{\bf Example.}\quad
For each $i\in I,$
let $\bone_i\in B=\k^I$ denote the idempotent corresponding to the $i$-th
direct
summand. Clearly, giving a finite rank $B$-bimodule amounts to
giving a finite dimensional $\k$-vector space $E$ equipped with
a direct sum decomposition: $E=\bigoplus_{i,j\in I}\;E_{i,j}\,,$
where $E_{i,j} = \bone_i\cdot E\cdot\bone_j$. 
Thus, one may think of the data $(B,E)$ as an oriented graph with vertex
set $I$ and with $\dim E_{i,j}$ edges going from the vertex $i$ 
to the vertex $j$.

Conversely, let $Q$ denote an oriented quiver with vertex set $I$.
Set $B=\k^I$, and let $E_Q$ be the $\k$-vector space
with basis formed by the set of edges $\,\{a\in Q\}\,.$ Then $E_Q$ has  an
obvious $B$-bimodule structure, and $T_{\!_B}(E_Q)$ is known as
the {\it path algebra} of $Q$. Further,
the $B$-bimodule $E_\dq$ associated with $\dq$, the {\it double} of $Q$,
has a natural symplectic $B$-form. The corresponding class in
$\dis\dr^2_{\!_B}\left(T_{\!_B}(E_Q)\right)$ is given by the
formula: $\;\om_{_\dr}=\sum_{a\in Q}\; da\otimes da^*\,.$\qed
\bigskip

In the special case $B=\k$, the  Lie bracket 
$\{\,,\,\}_{_{\om_{_\dr}}}$ on $A/[A,A]$
has been introduced by 
Kontsevich  [K2] in a somewhat different way as follows. 
Let $x_1,\ldots,x_n\,,\,y_1,\ldots,y_n$ be a symplectic basis of
the vector space $E$, i.e. a $\k$-basis such that:
$\om(x_i,y_j)=\delta_{ij},$
and $\om(x_i,x_j)= \om(y_i,y_j)=0$. 
By Lemma \ref{DR(quiver)}\vii, one has: $\dr^1_{\k}A \simeq A\otimes
E$. Kontsevich 
exploits this isomorphism
to write any 1-form $\alpha\in \dr^1_{\k}A$ in the form: 
$\alpha= \dis
\sum\nolimits_{i=1}^n F_{x_i}(\alpha)\otimes x_i + \sum\nolimits_{j=1}^n
 F_{y_j}(\alpha)\otimes y_j\,$, for certain uniquely determined elements
$ F_{x_i}(\alpha)\,,\,F_{y_j}(\alpha)\in A$. 
He then introduces , for any $i=1,\ldots n,\,$ the following $\k$-linear maps:
$$\frac{\partial}{\partial x_i}\,\,,\,\frac{\partial}{\partial y_i}:
\;\, \dr^0_{\k}A \too A\,,
\quad\mbox{\sl given by}\quad \frac{\partial f}{\partial {x_i}}:= F_{x_i}(df)\quad
\mbox{\sl and}\quad\frac{\partial f}{\partial {y_i}}:= F_{y_i}(df)
\,.
$$
Using these maps Kontsevich defines (put another way:
gives a coordinate expression for) the  Lie bracket $\{-, -\}\bom$
 by the familiar
formula:
\begin{equation} \label{poisson}
\{f, g\}_{_\om} \;:=\; \sum_{i=1}^n
\Bigl(\frac{\partial {f}}{\partial x_i}\cdot\frac{\partial {g}}{\partial y_i}
-\frac{\partial {f}}{\partial y_i}\cdot\frac{\partial {g}}{\partial x_i}\Bigr)
\;\mbox{\sf mod } [A,A]\enspace
\in\; A/[A,A] =\dr^0_{\k}A\,,
\end{equation} 
where "dot" stands for the product in $A$. We leave to the reader to 
check that  formulas (\ref{bom}), and
(\ref{mybracket}) give rise to the same bracket on $\dr^0_{\k}A=A/[A,A]$
as formula (\ref{poisson})
\bigskip

\noindent
{\bf Remarks.}\quad\vi In the general case of an arbitrary quiver
$Q$,  the analogue of Kontsevich's formula (\ref{poisson}) for the
Poisson bracket associated with the corresponding algebra $A=T_B(E_\dq)$,
 in  obvious notation, cf. (\ref{moment_quiver}), is:
\begin{equation} \label{quiv_K}
\{f, g\}_{_\om} \;=\; \sum\nolimits_{a\in Q}\;
\Bigl(\frac{\partial {f}}{\partial a}\cdot\frac{\partial {g}}{\partial a^*}
-\frac{\partial {f}}{\partial a^*}\cdot\frac{\partial {g}}{\partial
a}\Bigr)\;\mbox{\sf mod } [A,A]\enspace\in\; A/[A,A] =\dr^0_{_B}A\,.
\end{equation}

\vii The Poisson structure given by (\ref{quiv_K}) admits a
natural {\it quantization}, in which the symplectic
manifold $\Rep(\dq, V)=T^*\Rep(Q,V)$ gets replaced
by the algebra $\dd(\Rep(Q,V))$ of polynomial
differential operators
on the vector space $\Rep(Q,V)$. This quantization has been 
found by M. Holland [Ho], even before the non-commutative
Poisson structure given by (\ref{bom}) and (\ref{quiv_K})
has been discovered.\qed\bigskip

The next Proposition gives a non-commutative analogue 
of the  classical Lie algebra exact sequence:
$$
0\to\; \mbox{\it constant functions}\;  \to\;  \mbox{\it regular
functions} \; \to\; \mbox{\it 
symplectic vector fields}\;  \to\; 0\,,
$$
 associated
with a connected and simply-connected symplectic manifold.
\begin{proposition}\label{Lie}
\label{central_ext} There is a natural Lie algebra
central extension:
$$0\too B\too  A /[ A , A ]\too \derb( A ,\om)\too 0\,.$$
\end{proposition}

{\sl Proof.} It is immediate from formula (\ref{Psi}) 
that for the map: $f\mapsto \theta_f$ we have:
$\,\Ker\lbrace{A /[ A , A ]\too \derb( A ,\om)\rbrace}\,$
$=\,\Ker\,d$. By
Lemma \ref{DR(quiver)}\vi we get: $\Ker\,d=B$.
Further, Lemma \ref{DR(quiver)}\vi insures that
every closed element in $\dr^1_{\!_B}A$ is exact.
This yields surjectivity of the map:
$A /[ A , A ]\too \derb( A ,\om).$

It remains to show that the map: $f\mapsto \theta_f$
is a Lie algebra homomorphism. To this end,
we use the notation:
 $\delta= [\theta_f, \theta_g] - \theta_{\{f,g\}}
\in \derwa$
 from the proof of  Theorem \ref{Lie*}.
The   proof of  Theorem \ref{Lie*} implies
that the derivation $\delta$
has the property that,
for any $\eta\in \derb{A}$, one has:
$i_\theta i_\eta\om =0.$ Equivalently, setting $\al:= i_\delta\om$,
we get: $i_\eta\al=0\,,\,\forall \eta\in \derb{A}$.
But the isomorphism:
$\,\dr_{\!_B}^1(T_{\!_B}E) = (T_{\!_B}E)\bigotimes_B\;E\,$
of Lemma \ref{DR(quiver)}(ii) clearly forces such a 1-form
$\al$ to vanish.
The 2-form $\om$ being non-degenerate, it follows that
$\delta=0$. Thus, we have shown that
$\,[\theta_f, \theta_g] = \theta_{\{f,g\}}\,.$ 
This completes the proof.
\qed
\bigskip

\noindent
{\bf Representations.}\quad We now 
fix a finite dimensional left $B$-module $V$,
as at the end of \S2. Observe that if $Q$ is a quiver, and
$E=E_{\dq}$ is the symplectic $B$-bimodule attached,
as has been explained earlier, to the double of
$Q$,
 then for $A= T_{\!_B}(E_{\dq})$, 
in the notation of the Introduction we have:
$\Hom\balg(A, \End V)= \Rep(\dq, V)$.
In general, let $E$ be a finite dimensional symplectic  $B$-bimodule.
Then, for $A= T_{\!_B}E$, one has:
$\Hom\balg(A, \End V)=\Hom_{_{B{{\sf{\tdash}bimod}}}}(E,\End V)\,.$
The latter space can be naturally identified with
$\,E^\vee\otimes_{\!_B}\End V\,.$
Note that we have the symplectic $B$-form $\om^\vee$ on $E^\vee$,
and  a non-degenerate symmetric
bilinear form $\tr: \End V\otimes_{\!_B}\End V \too B,\,$
given by: $
\, (F_1, F_2) \mapsto \tr(F_1\ccirc F_2)\,.$
By standard Linear Algebra, the  tensor product  of a
 skew-symmetric and symmetric
 non-degenerate forms gives the skew-symmetric non-degenerate bilinear
form:
$\,\om_{_\rep}:=
\om^\vee\otimes\tr.\,$  The 2-form $\om_{_\rep}$
makes $\,E^\vee\otimes_{\!_B}\End V$,
hence, $\Hom\balg(A, \End V)$,
a symplectic  $B$-bimodule, therefore gives rise
to a $G(V)$-invariant Poisson bracket $\{\,,\,\}_{\om_{_\rep}}$ on the
coordinate ring $\k[\Hom\balg(A, \End V)]=\rep(A,V)$.
The invariants, ${\rep(A,V)}^{G(V)},$
clearly form a Poisson subalgebra in $\rep(A,V)$, and we have:

\begin{proposition}\label{omL}
The   map $\htr:  A /[ A , A ] \too {\rep(A,V)}^{G(V)}$
defined in (\ref{rep_map})
is a Lie algebra homomorphism, that is,
for any $f,g\in  A $, one has: 
$$\{\htr f\,,\,\htr g\}_{_{\om_{_\rep}}} =\htr(\{f,g\}_{\om_{_\dr}})\,.$$
\end{proposition}

{\sl Proof.} Straightforward calculation for $f,g$
taken to be non-commutative monomials.~$\square$
\medskip

We can now complete the proof of Theorem \ref{main_intro}.
As we have mentioned in the Introduction,
 the $G(V)$-action on $\Hom\balg(A, \End V)=
E^\vee\otimes_{\!_B}\End V\,$ turns out to be Hamiltonian,
and the corresponding moment map
$\mu: E^\vee\otimes_{\!_B}\End V\, \too \g(V)^* = \sg(V)$,
cf. (\ref{moment_quiver}),
is given by the following formula:
\begin{equation}\label{moment_quiver2}
\sum\nolimits_i\;\phi_i\otimes F_i\; \mapsto\; 
\sum\nolimits_{j<k}\;\om^\vee(\phi_j,\phi_k)\cdot [F_j,F_k]\;\in \;\sg(V)
\quad,\quad \phi_i\in E^\vee\,,\,F_i\in \End V\;.
\end{equation}
Fix $\O \subset \sg(V)$, a  closed $\ad\, G(V)$-orbit, and assume that the
group $G(V)$ acts freely on the subvariety $\mu^{-1}(\O) \subset
E^\vee\otimes_{\!_B}\End V$. Then, the orbit space $\mu^{-1}(\O)/G(V)$
is a smooth affine subvariety in $\Spec ({\rep(A,V)}^{G(V)})$.

\begin{proposition}\label{omL2}
The composite map:
$$\mu^{-1}(\O)/G(V)\;\into\; \Spec \bigl({\rep(A,V)}^{G(V)}\bigr)\;
\;\;\stackrel{{}_{\tiny{\sf evaluation}}}{\too}\;\;\;
\bigl({\rep(A,V)}^{G(V)}\bigr)^* \;\stackrel{\tr^*}{\too}\; ( A /[ A , A ])^*$$
 is injective and makes $\mu^{-1}(\O)/G(V)$ a coadjoint orbit in
$( A /[ A , A ])^*$.
\end{proposition}

{\sl Proof.} Set $X=\mu^{-1}(\O)/G(V)$, a smooth affine variety.
As we have argued in \S1, proving the proposition
amounts to showing that regular functions on
$\dis\Spec \bigl({\rep(A,V)}^{G(V)}\bigr)$ of the form $\tr(\hat{a})\,,\, a\in A/[A,A],\,$
separate points and tangents of the variety $\;X\; \subset\;$
$\dis \Spec
\bigl({\rep(A,V)}^{G(V)}\bigr)$.
This is clearly true for the whole algebra ${\rep(A,V)}^{G(V)}$,
since it is true for the  algebra $\k[X]$,
 and every regular function on $X$ is obtained from
an element of ${\rep(A,V)}^{G(V)}$, by restriction.

We now use the result of Le Bruyn- Procesi [LP], saying
that the algebra ${\rep(A,V)}^{G(V)}$ is generated by 
elements of the form:
$$\,\tr(\hat{\bone}_i\cdot\hat{x}_1\cdot\hat{x}_2\cdot\ldots\cdot\hat{x}_k
\cdot\hat{\bone}_i)\qquad,\qquad
i=1,\ldots,p\;,\;\; x_j\in E\;,\;\;k\geq 1\,.$$
The expression above is nothing but $\tr(\hat{a}),$ for
$\,a={\bone}_i\cdot({x}_1\otimes\ldots\otimes{x}_k)\cdot{\bone}_i
\in T^k_{\!_B}E\subset A.\,$
It follows that,
although the map
$\htr:  A /[ A , A ] \too {\rep(A,V)}^{G(V)}$
is not itself surjective,
the algebra ${\rep(A,V)}^{G(V)}$ is
generated by its image. Thus, elements of the image
separate points and tangents of the variety $X$.
\qed\bigskip

\section{Stabilization: infinite dimensional limit}
\setcounter{equation}{0}

We keep the setup of \S2; in particular, we let $B=\k^I$
and fix $A$,  a finitely generated associative $B$-algebra.
Any imbedding: $V\into V'$ of finite rank left $B$-modules
 induces a map:
$\Hom\balg(A, \End V)\into \Hom\balg(A, \End V')$, which is
a closed imbedding of affine algebraic varieties. The latter imbedding
gives
rise to the restriction homomorphism of coordinate rings
\begin{equation}\label{restriction}
r_{_{V',V}}:\;\;
\rep(A,V')^{G(V')}\too
\rep(A,V)^{G(V)}.
\end{equation}
Observe that the collection of
all finite rank  $B$-modules, $V$, forms a direct
system with respect to  $B$-module imbeddings, and we
set $\,V_{\!_\infty}:=
\underset{{^{\underset{V}{\too}}}}{\lim}\;V\,,$ and let
 $\,G_{\!_\infty}:= \underset{{^{\underset{V}{\too}}}}{\lim}\;G(V)\,$
be the corresponding  ind-group.
By definition we put:
 $\;\rep(A, V_{\!_\infty})^{G_{\!_\infty}} := 
\underset{{^{\underset{V}{\longleftarrow}}}}{\lim}\;\rep(A,V)^{G(V)}.\,$

There is a standard way to introduce
a cocomutative coproduct $\,
\Delta:\rep(A, V_{\!_\infty})^{G_{\!_\infty}} 
\too \rep(A, V_{\!_\infty})^{G_{\!_\infty}}\,\bigotimes_{_B}\,
\rep(A, V_{\!_\infty})^{G_{\!_\infty}}.\,$
To see this, it
is convenient to think of $\rep(A, V_{\!_\infty})$
  as some sort of coordinate ring
$\k[\Hom\balg(A,\End V_{\!_\infty})]$. Then any choice of $B$-module
isomorphism: $\,V_{\!_\infty} \stackrel{\varpi}{\simeq}
 V_{\!_\infty}\oplus V_{\!_\infty}\,$
gives a morphism of ind-schemes: 
$$\Hom\balg(A,\End V_{\!_\infty}) \times \Hom\balg(A,\End V_{\!_\infty})
\into \Hom\balg\bigl(A,\End(V_{\!_\infty}\oplus V_{\!_\infty})\bigr)
\stackrel{\varpi}{\stackrel{_\sim}{\to}} 
\Hom\balg(A,\End V_{\!_\infty}).$$ 
The coproduct $\Delta$ 
on $\rep(A, V_{\!_\infty})^{G_{\!_\infty}}$ is the one
induced by the corresponding algebra map:
$$\Delta:\;\;\k[\Hom\balg(A,\End V_{\!_\infty})]\;\; \too\;\;
\k[\Hom\balg(A,\End V_{\!_\infty})]\;\;\bigotimes\nolimits_B\;\;
\k[\Hom\balg(A,\End V_{\!_\infty})]\,.$$
Let $\prim\left(\rep(A, V_{\!_\infty})^{G_{\!_\infty}}\right)$
 denote the $B$-module of primitive elements in
  $\rep(A, V_{\!_\infty})^{G_{\!_\infty}}$,
i.e., the elements
$f\in \rep(A, V_{\!_\infty})^{G_{\!_\infty}}$
 such that $\,\Delta(f)= f\otimes 1 + 1\otimes f\,.$
\medskip

Observe further that the 
map $\,\htr_{_V} : A/[A,A]\too {\rep(A,V)}^{G(V)}$
given by (\ref{rep_map}) is compatible with 
restriction morphisms $r_{_{V',V}}$, see (\ref{restriction}), that is, for any
imbedding $\,V\into V'$, one has a commutative triangle:
$\,
\displaystyle r_{_{V',V}}\,\ccirc\, \htr_{_{V'}}
= \htr_{_V}.\,$
Therefore, the maps $\,\{\htr_{_V}\}\,$ give rise to a well-defined
  limit map 
$\,\dis\htr_{{_{\tiny\infty}}} : 
\A/[\A,\A] \too \rep(A, V_{\!_\infty})^{G_{\!_\infty}}\,.$

We now specialize to the setup of \S3 and assume
that $A=T_{_B}E$, for a certain finite rank projective
$B$-bimodule $E$.
The following result is, in
a sense, dual to the well-known relationship,
see [LQ], [L], between cyclic
homology of an associative
algebra $A$ and primitive homology of the Lie algebra $\glinfty(A)$.

\begin{proposition}\label{zelmanov}
 For  $A=T_{_B}E$, the map $\tr_{{_{\tiny\infty}}}$ 
sets up a bijection:
\[\tr_{{_{\tiny\infty}}}:
A/(B+[\A,\A]) \iso \prim(\rep(A, V_{\!_\infty})^{G_{\!_\infty}})\;.\]
\end{proposition}

 Notice next that, for $A=T_{_B}E$, we have:
 $\,\rep(A,V_{\!_\infty})=\k[\Hom_{_{B{{\sf{\tdash}bimod}}}}(E,\End V_{\!_\infty})],\,$
is a polynomial
algebra with a natural grading, that also induces a grading
on $\rep(A, V_{\!_\infty})^{G_{\!_\infty}}$.
Furthermore, the coproduct $\Delta$ is compatible
with the (graded) algebra structure, hence makes 
$\rep(A, V_{\!_\infty})^{G_{\!_\infty}}$ a commutative
and  cocomutative graded Hopf $B$-algebra. The structure theorem for commutative
and cocommutative graded
Hopf
algebras implies that $\rep(A, V_{\!_\infty})^{G_{\!_\infty}}$
 must be the symmetric algebra (over $B$)
on 
the  $B$-bimodule of   its primitive elements.
Therefore, Proposition \ref{zelmanov}  yields

\begin{corollary}\label{coproduct}
For  $A=T_{_B}E$, the map $\tr_{{_{\tiny\infty}}}$ extends, by multiplicativity,
to a graded   {\sf isomorphism} of Poisson algebras:
$\dis\,{\mathtt{Sym}}^{^{_\bullet}}\bigl(\A/(B+[\A,\A])\bigr)\;\; \iso\;\;
\rep(A, V_{\!_\infty})^{G_{\!_\infty}}\,.$\qed
\end{corollary}

\noindent
{\bf Remark.}\quad It is interesting to note that, for any
finite dimensional $V$ such that $\dim V > 1$, the variety
$\,\dis\Spec\bigl(\rep(A, V)^{G(V)}\bigr)\,$ is 
 quite
complicated, e.g., in the Calogero-Moser case.
  Nonetheless,  Corollary \ref{coproduct} says that
the `limiting' variety
$\,\dis\Spec\bigl(\rep(A, V_{\!_\infty})^{G_{\!_\infty}}\bigr)\,$
is always a vector space.
\bigskip

{\sl Proof of Proposition \ref{zelmanov}.} It is clear from definitions,
that $\tr_{_{\tiny\infty}}(f) \in \rep(A, V_{\!_\infty})^{G_{\!_\infty}}$
 is a primitive element, for any homogeneous element $f\in
\A$ such that $\deg f >0$.
 Furthermore, one 
verifies that any element not contained in the
image of the map $\tr_{{_{\tiny\infty}}}$
cannot satisfy the  equation
$\Delta(f)=f\otimes 1 + 1\otimes f$, hence, is not primitive. 
Thus, the map $\tr_{{_{\tiny\infty}}}$ is surjective, and it suffices to
prove it  is injective.

In order to avoid 
complicated notation, we restrict ourselves to
proving injectivity in the special case of the quiver $Q$ consisting
of a single vertex and a single edge-loop, that is
the Calogero-Moser quiver (the general case goes in a similar
fashion with minor modifications). Thus, we assume
that $A=\k\langle x,y\rangle$,
and therefore, $\,\rep(A, V_{\!_\infty})=
\k[\glinfty\oplus\glinfty],\,$ where
$\,\glinfty
:=\underset{\longrightarrow}{\lim}\;\,\gl_n(\k)\,.$
We must show that,
given $f\in A$, the equation: $\tr_{{_{\tiny\infty}}}(f)=0$
implies: $\,f\in [A,A]$. This is proved as follows
(the argument below
 seems to be standard, but we could not find an appropriate reference
in the literature).

Let ${\mathbb{A}}=\k\langle x_1, x_2,\ldots\,,\,y_1,y_2 \ldots\rangle$ 
be the free associative algebra on countably many variables,
and $[{\mathbb{A}},{\mathbb{A}}]$ the $\k$-linear 
subspace of ${\mathbb{A}}$ spanned by the commutators.
Similarly to formula (\ref{trace}), to any element $F\in {\mathbb{A}}$ one
assignes a polynomial function $\tr F$ in infinitely many 
matrix variables: $X_1,X_2,\ldots\,,\,Y_1,Y_2 \ldots\;\in\;\glinfty\,,$
by inserting matrices instead of formal variables.
We claim that: {\it
if $F$ is {\sl multi-linear} in all its variables,
and
the polynomial function $\tr F$ is identically
zero on $\glinfty$, then $F\in [{\mathbb{A}},{\mathbb{A}}]$.} To prove this,
note that modulo $[{\mathbb{A}},{\mathbb{A}}]$ we can write:
$\,F(x_1, x_2,\ldots\,,\,y_1,y_2 \ldots)=
x_1\cdot \widetilde{F}(x_2,\ldots\,,\,y_1,y_2 \ldots)\,.$
Hence, equation: $\dis 0=\tr F(X_1, X_2,\ldots\,,\,Y_1,Y_2 \ldots)
=\tr\bigl(X_1\cdot \widetilde{F}(X_2,\ldots\,,\,Y_1, Y_2 \ldots)\bigr)$
implies, since the trace pairing on $\glinfty$
is nondegenerate, 
 that the function $\widetilde{F}(x_2,\ldots\,,\,y_1,y_2 \ldots)$
is identically zero on $\glinfty$. Furthermore, since $\glinfty$
(viewed as an associative algebra) is known to be
an algebra without polynomial identities, 
we conclude that $\widetilde{F}=0$. Thus,
$F\in [{\mathbb{A}},{\mathbb{A}}]$, and our claim is proved.

We can now complete the proof of the Proposition.  Fix
$f\in\A$ such that $\tr f(X,Y)=0$ identically on $\glinfty\oplus\glinfty$.
Rescaling transformations: $X\mapsto t\cdot
X\,,\, Y\mapsto s\cdot Y\;,\,\forall
t,s\in \k^\times,\,$ show that we may reduce to the
case where $f$ is homogeneous in $X$ and $Y$ of degrees, say
$p,q$, respectively. We now use the standard polarisation trick, and
formally
substitute: $x= t_1\cdot x_1+\ldots t_p\cdot x_p\;,\;
y= s_1\cdot y_1+\ldots s_q\cdot y_q\,$ into $f$,
and then take the term multilinear in $t_1,\ldots,s_q$.
This way we get from $f\in \A$ a multilinear element
$F\in {\mathbb{A}}$ such that $\tr F =0$ identically on $\glinfty$.
By the claim of the preceeding paragraph we conclude
that $F\in [{\mathbb{A}},{\mathbb{A}}]$.
Observe now that sending all the $x_i$'s to $x$, and 
all the $y_i$'s to $y$ yields an
algebra homomorphism $\pi: {\mathbb{A}} \to \A$ such that
$\pi(F)=p!q!\cdot f$. Applying this homomorphism to $F$ we get:
$f=\frac{1}{p!q!}\pi(F) \in \pi([{\mathbb{A}},{\mathbb{A}}]) =
[\A,\A]\,.$\qed\bigskip


\section{The basics of $\PP$-geometry.}
\setcounter{equation}{0}

 Let 
$\PP=\{\PP(n),\,n=1,2,\ldots\}$ be a $\k$-linear quadratic operad
with $\PP(1)=\k$, see
[GiK]. Let $\s_n$ denote the Symmetric group on $n$ letters.
Given $\mu\in \PP(n)$ and a $\PP$-algebra $A$,
we will write: $\mu_A(a_1,\ldots,a_n)\,$
for the image of $\mu\otimes a_1\otimes\ldots\otimes a_n$
under
the structure map: $\PP(n)\otimes_{_{\s_n}} A^{\otimes n} \too A\,.$
Following [GiK, \S1.6.4], we introduce
an {\it enveloping algebra} $\upa$, the associative unital
$\k$-algebra  such that
the abelian category of (left) $A$-modules is equivalent to the
category of left modules over  $\upa$, see [GiK, Thm. 1.6.6]. 
The algebra $\upa$ is generated by the symbols:
$u(\mu, a)\,,\, \mu\in \PP(2), a\in A,$ subject to certain
relations, see [Ba, \S1.7].

An ideal $I$ in a $\PP$-algebra $A$ will be called $N$-nilpotent if,
for any $n\geq N\,,\,\mu\in\PP(n)$, and $a_1,\ldots,a_n\in A\,,$ one
has: $\mu_A(a_1,\ldots,a_n)=0,$ whenever at least $N$ among the elements
$a_1,\ldots,a_n$ belong to $I$.
The following useful reformulation of the notion
of a left $A$-module is essentially well-known, see e.g., [Ba,~1.2]:

\begin{lemma}\label{La-mod}
Giving a left $A$-module structure on a vector space $M$ is equivalent
to giving a $\PP$-algebra structure on $\,A\sharp M:= A\oplus M\,$ 
such that the  following
conditions hold:

\vi The imbedding: $a\mapsto a\oplus 0$
makes $A$ a $\PP$-subalgebra in $A\sharp M$.

\vii $\;M$ is a 2-nilpotent ideal in $A\sharp M$. \quad\qed
\end{lemma}

A $\PP$-algebra in the monoidal category of
$\Z/2$-graded, (resp. $\Z$-graded) super-vector spaces, 
see [GiK, \S1.3.17-1.3.18],
will be referred to as
a  $\PP$-{\it superalgebra}, (resp. graded
 superalgebra). Any $\PP$-algebra may be regarded as a
$\PP$-superalgebra concentrated in degree zero.
Given a finite dimensional (super-) vector space $V$,
write $\overline{V}$ for the same vector space with reversed parity.
Let
$$
\T^{^{_\bullet}}\pp\!{V} :=\bigoplus\nolimits_{i\geq 1}\;\;
\PP(i) \otimes_{_{\s_i}}V^{\otimes i}\;\aand\;
\Tc^{^{_\bullet}}\pp{_{\!}}{V} :=\bigoplus\nolimits_{i\geq 1}\;\;
\PP(i) \otimes_{_{\s_i}}\overline{V}^{\otimes i}
\;$$
be the free graded $\PP$-algebra 
(resp. super-algebra) generated by $V$.

Fix a $\PP$-algebra $A$, and consider the category
of $A$-algebras, i.e.
of pairs $(B,p)$, where $B$ is a $\PP$-algebra and
$p: A\to B$ is a $\PP$-algebra morphism.
Note that such a  morphism makes $B$ an $A$-module.
Thus, we get an obvious forgetful functor:
$A${\it-algebras} $\too$ $A${\it-modules}.
The result below says that this functor
has a right adjoint:

\begin{lemma}\label{TM}\vi  Given a $\PP$-algebra $A$, there
is a functor: $M \mapsto T_A^{^{_\bullet}}M\,,$
(resp. $M \mapsto \tc^{^{_\bullet}}_AM$) assigning to a left
$A$-module $M$ a graded 
 $\PP$-algebra 
$\,T^{^{_\bullet}}_AM= \bigoplus_{i\geq 0}\;T^i_AM$ (resp. graded 
 $\PP$-superalgebra $\,\tc^{^{_\bullet}}_AM
= \bigoplus_{i\geq 0}\;T^i_A\overline{M}\,$), such that
$T^0_AM=A$.

\vii For any $\PP$-algebra map: $A\to B$, one has
 a natural adjunction isomorphism:
$$\Hom_{_{A{\sf{\tdash}mod}}}(M,B)\iso 
\Hom_{_{\PP{\sf{\tdash}alg}}}(T^{^{_\bullet}}_AM, B)\,.$$
\end{lemma}

{\textsl  {Proof}:} If $A$ is a $\PP$-subalgebra in a 
$\PP$-algebra $\widetilde{A}$,
we define a $\PP$-algebra $T_{_A}\widetilde{A}$ as the quotient of
$\T\bulp{\widetilde{A}}$, a free $\PP$-algebra, modulo two-sided ideal
generated by all relations of the form:
$$\mu\otimes a\otimes \tilde{a} = \mu_{_{\widetilde{A}}}(a, \tilde{a})\quad,
\quad \mu\otimes \tilde{a}\otimes a =
 \mu_{_{\widetilde{A}}}(\tilde{a},a)\quad,\quad \forall \mu \in
\PP(2),\,
a\in A\sset {\widetilde{A}},\, \tilde{a}\in {\widetilde{A}}\,,$$
where 
$\,\mu\otimes \tilde{a}\otimes a\,,\,\mu\otimes a\otimes \tilde{a} \in
 \PP(2)\otimes {\widetilde{A}}^{\otimes 2} = 
\T\pp^2{\widetilde{A}}$, and 
$\, \mu_{_{\widetilde{A}}}(a, \tilde{a}),\,\mu_{_{\widetilde{A}}}(\tilde{a},a)
\in
 \PP(1)\otimes {\widetilde{A}} = \T\pp^1{\widetilde{A}}$.
 We now apply this construction to
the algebra ${\widetilde{A}}=A\sharp M$,
and put $\,T^{^{_\bullet}}_AM := T_{_A}{\widetilde{A}}$,
where the grading on the left accounts for the number of occurrences
of elements of $M$, which is well-defined since the
relations involved in the definition of $T_{_A}{\widetilde{A}}$ are `homogeneous in
$M$'. 

A closer look at the construction above shows that
\begin{equation}\label{T_AM}
T_A^{^{_\bullet}}M \;=\;A\;\,\bigoplus\;\,
(\T\bulp{M})\big/\langle\!\langle\mu^{(12)}(a, m_1) \,\otimes\,m_2 -
m_1\,\otimes\mu(a,m_2)\rangle\!\rangle\,,
\end{equation}
where $\langle\!\langle\ldots\rangle\!\rangle$ denotes
the two-sided ideal generated by the indicated
subset of $\PP(2)\otimes M^{\otimes 2}=
\T^2\pp{M}$, for all
$\mu\in\PP(2), a\in A, m_1,m_2\in M\,,$ and where
$\mu^{(12)}$ stands for the action of the transposition
$(12)\in\s_2$
on $\mu$.
In particular, we have:
 $T^0_AM=A$ and $T^1_AM=M$.\qed\medskip

Let $A$ be a $\PP$-algebra and $M$ a left $A$-module. By Lemma 
\ref{La-mod}, we may (and will) regard 
$A\sharp M$ as a
$\PP$-algebra.

\begin{definition}\label{der}
 A $\k$-linear map $\theta:
A\to M$ is called a {\it derivation} if
the map: $\,a\bigoplus m \mapsto a\,\bigoplus\, \theta(a)\!+\!m\,,$
is an automorphism of the $\PP$-algebra $A\sharp M$.
\end{definition}
\noindent
Equivalently, following [Ba, Definition 3.2.6],
 extend $\theta$ to a $\k$-linear map $\theta^\sharp:
A\sharp M \to A\sharp M\,,$ given by 
$\,\theta^\sharp: a\oplus m \;\mapsto \;0\oplus \theta a$.
Then, $\theta$ is a derivation if and only if,
for any  $\mu\in \PP(n)$, we have:
$$\theta^\sharp\bigl(\mu_{_{A\sharp M}}(b_1,\ldots,b_n)\bigr)
=\sum_{i=1}^n\; \mu_{_{A\sharp M}}(b_1,\,\ldots,\, b_{i-1},\, \theta^\sharp{b_i}\,,\,
b_{i+1},\,\ldots,\, b_n)\quad,\enspace\forall
b_1,\,\ldots,\, b_n\in A\sharp M.$$

Let
$\derp(A,M)$ denote the $\k$-vector space of all derivations
from $A$ to $M$. It is straightforward to see that 
 the ordinary commutator makes $\derp(A,A)$ 
a Lie algebra. \medskip

Next we define, following [Ba, Def. 4.5.2], 
an $A$-module of {\it K\"ahler differentials}
 as the left $\upa$-module, $\Om^1\pp{A}$, generated by the symbols
$da$, for $a\in A$, subject to the relations:

\vi $\quad\;\; d(\lambda_1 a_1 +\lambda_2 a_2) = \lambda_1 da_1+\lambda_2
da_2\quad,\enspace\forall \lambda_1,\lambda_2 \in \k;$

\vii $\quad
 d(\mu(a_1,a_2)) = u(\mu,a_1)\otimes da_2 + u(\mu^{(12)}, a_2)\otimes
da_1\quad,\enspace
\forall \mu\in \PP(2)\,,\,a_1,a_2\in A\,,$ 

\noindent
where $u(\mu,a)$
denote the standard generators of $\upa$, see [Ba].
\medskip

By construction, $\Om\pp^1\!A$ is a left $A$-module,
and the assignment $a\mapsto da$ gives a derivation
$d\in \derp(A\,,\,\Om\pp^1\!A)$.
Moreover, this derivation is universal in the following
sense.
Given any left $A$-module $M$ and a derivation
$\theta: A\to M$, there exists an 
$A$-module morphism  $\,\Om^1\theta: \Om\pp^1\!A \to M$,
uniquely determined by the condition that 
$\,(\Om^1\theta)(da)=\theta(a)\,.$ It follows 
that the 
$A$-module of K\"ahler differentials
represents the functor
$\derp(A,-)$, i.e., we have (see [Ba, Remark 4.5.4]):

\begin{lemma}\label{om-der} 
For any left $A$-module $M$ there is a natural isomorphism:
$$\derp(A,M) \;\simeq\; \Hom_{_{A{\sf{\tdash}mod}}}(\Om^1\pp{A}, M)\,.\quad\square$$
\end{lemma}

In particular, for $M=A$, we get an isomorphism: $\derp(A,A)
 \iso \Hom_{_{A{\sf{\tdash}mod}}}(\Om^1\pp{A}, A)\,.$
We let $i_\theta\in\Hom_{_{A{\sf{\tdash}mod}}}(\Om^1\pp{A}, A)$ denote
the morphism: $\Om^1\pp{A}\to A$,
corresponding to $\theta\in \derp(A,A)$ under the  isomorphism above.

We set $\,\Om\bulp A := \tc^{^{_\bullet}}_A(\Om^1\pp{A}),\,$
a graded $\PP$-superalgebra generated by the $A$-module $\Om^1\pp{A}$.
\medskip
Recall that the {\it differential envelope}
of a $\PP$-algebra $A$ is a differential graded  $\PP$-super-algebra
$D^\bullet(A)= \bigoplus_{i\geq 0}\;D^i(A),\,$ such that
$D^0(A)=A$, and such that the following universal property holds:
{\it For any differential graded $\PP$-superalgebra
$\widetilde{D}^\bullet= \bigoplus_{i\geq 0}\;\widetilde{D}^i\,,$ and a $\PP$-algebra 
morphism $\rho: A\to \widetilde{D}^0$ there exists a unique
DG-superalgebra morphism $D(\rho): D^\bullet(A) \to  \widetilde{D}^\bullet$
such that $D(\rho)\big|_{_{D^0(A)}} = \rho$.}

\begin{proposition}\label{d^2=0} \vi On $\Om\bulp A$, there exists
 a natural
super-differential $\,d: \Om\bulp A \too \Om^{\bullet+1}\pp{A}\,,\,
d^2=0,\,$
such that its restriction: $A=\Om^0\pp{A}\too \Om^1\pp{A}$
coincides with the canonical $A$-module
derivation $d: A \to \Om^1\pp{A}$. 

\vii The differential graded $\PP$-superalgebra $\,(\Om\bulp A\,,\,d)\,$
is the differential envelope of $A$.
\end{proposition}

\noindent
{\textsl {Proof.}} We first give a direct construction of 
the differential envelope $D^\bullet(A)$ of a $\PP$-algebra $A$,
 as follows. Let $\ba$ denote a second copy of $A$
viewed as a $\k$-vector space,
and write $\overline{a}$ for the element
of $\ba$ corresponding to  an element $a\in A$.
We form the graded super-vector space $A\oplus \ba$, where
$A$ is placed in grade degree zero, and
$\ba$ is  placed in grade degree 1. 
Let
$\Tc^{^{_\bullet}}\pp{_{\!}}(A\oplus \ba) :=\bigoplus\nolimits_{i\geq 1}\;\;
\PP(i) \otimes_{_{\s_i}}(A\oplus \ba)^{\otimes i}
\;$ be the free $\PP$-superalgebra generated by $A\oplus \ba$,
viewed as a graded superalgebra with respect to 
the total grading coming from both
the grading on  $A\oplus \ba$ and the grading on the tensor algebra.
 We put: $D^\bullet(A) \;:=\;\Tc^{^{_\bullet}}\pp{_{\!}}(A\oplus
\ba)/I\,,$ where $I$ is the two-sided ideal
generated by the following set:
\begin{equation}
\label{dif_envelope}
\{\mu\otimes a_1\otimes a_2 - \mu(a_1,a_2)\;\;,\;\;
\mu\otimes\overline{a}_1\otimes a_2 +\mu\otimes a_1\otimes 
\overline{a}_2
-\overline{\mu(a_1,a_2)}\}_{\,\mu\in \PP(2)\,,\,
a_1,a_2\in A}.
\end{equation}
Thus, $D^\bullet(A)$ is a graded $\PP$-superalgebra.

The $\k$-linear endomorphism of  $A\oplus \ba$ given
by the assignment: $a\oplus \overline{a}_1 \mapsto 0\oplus \overline{a}$
extends uniquely to a super-derivation
$\Tc^{^{_\bullet}}\pp{_{\!}}(A\oplus \ba) \too
\Tc^{^{_\bullet}}\pp{_{\!}}(A\oplus \ba).$
This derivation descends to
a well-defined derivation $d$ on $D^\bullet(A)$. Note that,
for any $x\in A\oplus \ba$ we have: $d^2(x)=0$.
This implies, since the subspace $A\oplus \ba$ generates
the algebra  $D^\bullet(A)$, that  $d^2=0$ identically on $D^\bullet(A)$.
 Thus, $d$ makes  $D^\bullet(A)$  a differential graded
$\PP$- superalgebra. 

The zero-degree component, $D^0(A)$, of the super-algebra $D^\bullet(A)$
is by construction
a $\PP$-subalgebra isomorphic to $A$, i.e.,
there is a canonical superalgebra imbedding
$j: A=D^0(A) \into D^\bullet(A)$. Hence,
$D^\bullet(A)$ may be regarded as an $A$-module, and the assignment: 
$a\mapsto \overline{a}$ gives a derivation $d \in \derp(A, D^\bullet(A))$.
This derivation is universal in the sense explained above
(for uniqueness property use that the superalgebra $D^\bullet(A)$
is generated by the subspace
 $A\oplus \ba$).
Hence $D^\bullet(A)$ is the differential envelope of $A$.

Observe next that the degree 1 component
of $D^\bullet(A)$ is  isomorphic, by definition of $D^\bullet(A)$,
to the quotient of $\,\upa\otimes
\ba\,$ by the relations \vi\!--\vii defining the module $\Om\pp^1A$ of 
K\"ahler differentials. Therefore, $D^1(A)$, the degree 1 component of $D^\bullet(A)$,
is isomorphic to $\Om\pp^1A$ and, moreover, 
the canonical
derivation $d: A\to \Om\pp^1A$  may be identified with the
map: $a\mapsto \overline{a}\in D^\bullet(A)$.

By the universal property of the tensor algebra,
the $A$-module imbedding $\Om^1A \into D^\bullet(A)$ can be extended uniquely
to a graded super-algebra morphism
$f: \tc^{^{_\bullet}}_A(\Om^1\pp{A})\too D^\bullet(A)$.
To show that $f$ is an isomorphism we construct its
inverse, a map $g: D^\bullet(A)\too\tc^{^{_\bullet}}_A(\Om^1\pp{A})$,
as follows. We have an obvious imbedding of $\k$-vector spaces:
$\,A\oplus \ba \into A \oplus\Om^1\pp{A}$, given by:
$a\oplus \overline{a}_1 \mapsto a\oplus da_1\,.$ This imbedding
extends, by the universal property of a free $\PP$-algebra,
to a $\PP$-superalgebra morphism
$\,\tilde{g}: \Tc^{^{_\bullet}}\pp{_{\!}}(A\oplus \ba) \too
\tc^{^{_\bullet}}_A(\Om^1\pp{A}).$ The relations defining the ideal $I$
 in formula (\ref{dif_envelope}) are designed in such a way
that the morphism $\tilde{g}$ descends to a well-defined
super-algebra morphism $\,g: D^\bullet(A) \to
\tc^{^{_\bullet}}_A(\Om^1\pp{A}).$
It is straightforward to verify that $g=f^{-1}$.
\qed\medskip

\begin{remark} Our construction agrees with the notion of
non-commutative differential forms for an algebra over the associative
operad, as defined e.g. in [L] and  used in \S2 above.
\end{remark}
\bigskip

From now on we assume, in addition,  that
$\PP$ is
a  {\it cyclic} Koszul operad, see [GeK], with $\PP(1)=\k$.
In particular, for each $n\geq 1$, the space $\PP(n)$ is
equipped with an $\s_{n+1}$-action that extends
the $\s_n$-module structure on $\PP(n)$ arising from
the operad structure. Write ${\mathtt {Sym}}_{_\k}^2A$ for the symmetric square
of $A$.
Following an idea of Kontsevich, Getzler and Kapranov
introduce a functor
$\,\,\R: \mbox{\it $\PP$-algebras}\;\too\;\k\mbox{\it-vector spaces}\;,$
$$\R:
\; A\;\;\; \mapsto\;\;\; \R(A) := {\frac{{\mathtt  {Sym}}_{_\k}^2A}
{\big\langle a_0\cdot\mu(a_1,a_2)- \mu(a_0,a_1)\cdot
a_2\big\rangle}}_{\big|\; a_0,a_1,a_2 \in A\,,\,\mu\in\PP(2)}
\,.$$

Generalizing  the Karoubi's construction [Ka] in the associative case,
 define de Rham complex of $A$ as the 
graded vector space
$\;\dr^{^{_\bullet}}{A} \,:= \,\R(\Om\bulp A)\,.$ 
The differential $d$ on $\Om\bulp A$ induces a differential
on $\dr^{^{_\bullet}}{A}$.

For any $\theta \in \derp{A} $, the  morphism $i_\theta: \Om^1\pp{A}\to
A$ 
introduced after Lemma \ref{om-der}  extends 
 to a   super-derivation
$i_\theta:
\Om\bulp A\to\Om\pp^{\bullet -1\!}A$,
 called the  {\it contraction} operator.
Further, the derivation $\theta$  induces, by a standard argument,
 a derivation $L_\theta$ of the
associative algebra $\upa$, and a map
$\,L_\theta: \Om^1\pp{A} \to \Om^1\pp{A}$. The latter one extends
to a derivation $L_\theta:
\Om\bulp A\to\Om\bulp A$,
 called the {\it Lie operator}.
The maps $i_\theta$ and $L_\theta$ descend naturally to the corresponding operators
on $\dr^{^{_\bullet}}{A}$.
It is straightforward to verify that these latter operators 
satisfy all the  standard commutation relations (\ref{identities}). 

\bigskip

\section{Symplectic geometry of a free ${\mathbf{\mathscr{P}}}$-algebra.}
\setcounter{equation}{0}
We keep the assumption that $\PP$ is a cyclic Koszul operad.
In this section which is a generalization of \S3,
inspired by  works of Drinfeld [Dr, Proposition 6.1 and above it] and
Kontsevich (private communication, 1994), we consider 
the case of a free $\PP$-algebra. To avoid unnecessary
repetitions and to simplify notation we only consider the `absolute'
case, i.e., the case of the ground ring $B=\k$.

Fix a finite dimensional $\k$-vector space $E$, and write
$\A=\T\pp\!{E}$ for the free $\PP$-algebra (note that 
$\PP$-algebras are algebras without unit, in general).
We have:
\begin{equation}\label{iso}
\R(\A) = \dr^0(\A) =
\bigoplus\nolimits_{i\geq 1}\;\PP_i\otimes_{_{\s_{i+1}}}E^{\otimes(i+1)}
\quad,\quad
\dr^{1}(\A)\; =\; \A\; \bigotimes\;E\;.
\end{equation}

Let $\widehat{\A}= \prod_{i\geq 0}\, \T^i\pp\!{E}$ denote the completion of $\A$ with
respect to the augmentation, and let $\aut(\widehat{\A})$
denote the group of continuous  algebra automorphisms of $\widehat{\A}$.
Any such automorphism $\Phi$ is determined by its restriction to
$E= \T^1\pp\!{E}$, a $\k$-linear map $\phi: E\to \widehat{\A}$.
We have an expansion: $\phi(v)= \sum_{i=1}^\infty\phi_i(v),$ where
$ \phi_i(v) \in\T^i\pp\!{E}$. We write $d\Phi: E\to E$,
for the map: $v\mapsto \phi_1(v)$; and we let
$\aut_\circ(\widehat{\A})$ be the subgroup of
$\aut(\widehat{\A})$ formed by all automorphisms $\Phi$ such that
$d\Phi=\id_E$.

Observe further that the obvious grading on the free algebra $\A=\T\pp\!{E}$
induces a  natural
grading $\R^{^{_\bullet}}(\A)= \bigoplus_i\,\R(\A)_{(i)}$,
and, for each $p\geq 0$,
 a similar grading $\dr^p(\A)= \bigoplus_i\,\dr^p(\A)_{(i)}$.
Fix a closed 2-form $\om \in \dr^2\A$, and let
$\om =\om_0 + \om_1+\ldots,\,
\om_i\in \dr^2(\A)_{(i)}\,,$ be its expansion into graded
components. We see, in particular, that $\om_0$ may be viewed as an
ordinary
skew-symmetric $\k$-bilinear form: $E\times E \to \k$.

\begin{theorem}\label{Darboux}{\bf(Darboux theorem)}
$\;$
\vi A closed 2-form  $\om =\om_0 + \om_1+\ldots\in \dr^2\A$ is 
non-degenerate if and only if so is the associated bilinear form
$\om_0: E\times E \to \k$.

\vii If $\om$ is non-degenerate then there exists an automorphism
$\Phi\in \aut_\circ(\widehat{\A})$ such that:
$\Phi^*\om = \om_0.$
\end{theorem}

{\textsl {Proof.}} Part \vi is clear. Part \vii is proved by the
standard `homotopy argument'. Specifically, 
we consider a 1-parameter `family':$\,\om_t= \om_0 + 
t\cdot\om'\in \dr^2\A[[t]],$
where $\om'=\om-\om_0=\om_1+\om_2+\ldots\in \dr^2\A$. The 2-form $\om'$ being
closed, there exists $\alpha\in \bigoplus_{p\geq 1}
\,\dr^1(\A)_{(p)},$ such that $\om'= -d\alpha$.
Since $\om_0$ is non-degenerate, there exists a 1-parameter family 
$\,\theta_t \in \k[[t]]\,\widehat{\otimes}\,\derp\A=\derp\A[[t]]\,$ determined uniquely
from the equation: $\,i_{\theta_t}\om_t= \alpha$. We define
$\,\Phi(t)\in \aut(\widehat{\A}[[t]]),\,$
a formal one-parameter family of automorphisms of $\A$,
to be the solution of the differential equation:
$\,\frac{d\Phi(t)}{dt}=L_{\theta_t}\Phi(t)\,$ of the form:
$\,\Phi(t)=
\id_{_{\A}}+ t\cdot \Phi_1 + t^2\cdot \Phi_2 + \ldots\,.$
It follows from the construction that $\Phi(t)^*\om_t=\om_0$,
see  e.g. [GS] for more details.
Note further that
the series $\Phi(t)$ above has only finitely many terms in any given  grade
degree $p\geq 0$,
i.e. terms that shift the grading on $\A$ by $p$. In particular, setting $t=1$
in this series gives a well-defined element of $\Phi(1)\in
\aut_\circ(\widehat{\A})$
and we get: $\Phi(1)^*\om_{_{t=1}}=\om_0$.
But $\om_{_{t=1}}=\om$, and part (ii) follows.\qed\medskip

Because of this result, there is no loss of generality in 
considering only degree zero symplectic 2-forms $\om\in \dr^2\A$,
i.e., such that  $\om = \om_0$. Fix such an $\om$,  that is fix
$(E,\om)$, a symplectic vector space. Imitating  the strategy used
in \S2 it is possible to define a Lie bracket on the vector space
$\R(\A)$. We prefer however to give the following
 direct explicit construction of this bracket
similar to formula (\ref{mybracket}) in the associative case.

 For each $i,j \geq 1$, let
$\star: \mu \otimes \nu \mapsto \mu(1,\ldots,1,\nu)\,,$
denote the operad-composition map:
$\,\PP(i) \otimes \PP(j)\, \simeq\,\PP(1)\otimes \ldots\otimes\PP(1)
\otimes\PP(i) \otimes \PP(j) \too \PP(i+j-1)\,,$
where $1\in \k=\PP(1)$, see [GeK, Theorem 2.2(2)].
We now change the notation and write:
$\R^\bullet(A)=\bigoplus_i\,\R^i(A),$ where
$\R^i(A)$, previously denoted by
$\R(A)_{(i)}$, is the graded component
with respect to the grading induced by one on $A$.
Also, let ${\sf {Sym}}$ be the
`symmetrisation map', the projection to $\s_n$-coinvariants.
For each $i,j \geq 1$, we define
a bilinear pairing $\,\{-,-\}\bom:\,\R^i(\A) \otimes \R^j(\A)\too
\R^{i+j-1}(\A)\,$ as the following composition
\begin{align}\label{P-bracket}
 \R^i(\A) \otimes \R^j(\A) = &
\left( \PP(i)\otimes_{\s_{i+1}}E^{\otimes i+1}\right)\;\bigotimes\;
\left(\PP(j)\otimes_{\s_{j+1}}E^{\otimes j+1}\right) \too
\nonumber\\
&
\left(\PP(i)\otimes \PP(j)\otimes E^{\otimes i+ j+2}\right)_{\s_{i+1}\times \s_{j+1}}
\stackrel{\star}{\too}\nonumber\\
& \left(\PP(i+j-1)\otimes E^{\otimes i+ j+2}\right)_{\s_{i+1}\times
\s_{j+1}}
\stackrel{{\sf {Sym}}}{\;\too\;}\nonumber\\
&
\left(\PP(i+j-1)\otimes_{\s_{i+j}}
E^{\otimes i+ j}\right)\;\bigotimes\; E^{\otimes 2}\,
\stackrel{{\mathsf{id}}\otimes\om}{\too}\nonumber\\
& \PP(i+j-1)\otimes_{_{\s_{i+j}}}\; E^{\otimes i+ j}
\;=\;\R^{i+j-1}(\A)\,.\nonumber
\end{align}

An appropriate modification of the proof of Theorem \ref{Lie*}, or  a
 direct calculation, yields

\begin{proposition}\label{kontsevich_bracket} The bracket $\{-,-\}$
makes $\R^{\bullet-1}(\A)$ into a graded Lie algebra.\hfill\qed
\end{proposition}

Let $\derp(\A,\om)$ denote the Lie subalgebra in $\derp\A$ formed by all
derivations $\theta \in \derp\A$ such that $L_\theta\om=0$.
Since $\om=\om_0$, this is equivalent to the requirement
that the degree zero component  $d\theta: \A_1\to \A_1$ induces an
endomorphism of $E\otimes E$ that
annihilates $\om^\vee\in E\otimes E$. Using the same argument as in
\S\S2-3, one proves the following two results

\begin{lemma}\label{forms_fields} The assignment: $\theta \mapsto i_\theta\om$ gives
graded vector space isomorphisms:
$$ \derp^{^{_\bullet}}(\A) \iso \dr^1(\A^{^{_\bullet}})\aand
\derpa \iso \dr^1(\A^{^{_\bullet}})_{_{\sf{closed}}}\,.\eqno\square$$
\end{lemma}

\begin{proposition}\label{central_extension}
 There is a canonical  graded Lie algebra central extension:
$$0\too \k \too \R^{\bullet-1}(\A) \too \derpa\too 0\,.\eqno\square$$
\end{proposition}
\bigskip

We call a pair $(\sss, \tr)$, where $\sss$ is a $\PP$-algebra 
and $\tr$ is a symmetric non-degenerate invariant bilinear
form $\tr: \sss\otimes \sss
\to \k$, a {\it symmetric} $\PP$-algebra. 
 Any such bilinear form is determined, cf. [GeK], by a linear function
$\tr: \R(\sss)\to \k\,,\, b\otimes b'\mapsto \tr(b\otimes b')=\tr(b, b')$. 

From now on, fix
a finite-dimensional  symmetric $\PP$-algebra
$(\sss, \tr)$. Let $\aut(\sss, \tr)$ denote the algebraic group of
automorphisms of the $\PP$-algebra $\sss$ that preserve the
bilinear form $\tr$. The corresponding Lie
algebra $\derp(\sss, \tr)$ is formed by all the derivations
$\theta\in \derp{\sss}$ such that, for any $b,b'\in \sss$, one has:
$\,\tr(\theta(b)\,,\,b') +\tr(b\,,\,\theta(b'))=0\,.$
\medskip

\noindent 
{\textbf {Representation functor.}}\quad
For any finitely generated $\PP$-algebra $A$,
the set $\Hom\palg(A, \sss)$ has the natural structure of a finite
dimensional
affine algebraic variety, acted on by the algebraic group $\aut(\sss,
\tr)$. We put $\,\rep(A,\sss):= \k[\Hom\palg(A, \sss)].$

Let now $(E,\om)$ be a finite dimensional symplectic
vector space, and
$\A= T\pp{E}$, the free $\PP$-algebra on $E$.
Then we clearly have:
$\Hom\palg(\A, \sss) = \Hom_{\,_\k}(E, \sss) = E^*\otimes_{_\k} \sss$,
is a finite dimensional $\k$-vector space. 
The symplectic 2-form
$\om$ on $E$ gives rise, as in \S3, to the
symplectic 2-form $\,\om_{_\rep}:=\om^\vee \otimes\tr\,$ on 
$\Hom_{\,_\k}(E, \sss) = E^*\otimes_{_\k} \sss$. The action
of the group $\aut(\sss,\tr)$ on  $\Hom\palg(A, \sss)$ preserves this symplectic
form and is, moreover, Hamiltonian. 
In other words, the vector field
on $E^*\otimes_{_\k} \sss$  arising from a 
derivation $\theta\in \derp(\sss,\tr)$
is induced by an $\aut(\sss,\tr)$-invariant
Hamiltonian function $H_{_\theta} \in \rep(A,\sss)$.
Explicitly,
the function $H_{_\theta}$ is given by the following  quadratic polynomial 
on $E^*\otimes_{_\k} \sss$:
$$ H_{_\theta}:\;\;\sum\nolimits_k\;\, \check{x}_k\otimes s_k\;\;\mapsto
\;\;\sum\nolimits_{i< j}\;\, \om^\vee(\check{x}_i\,,\,\check{x}_j)\cdot
\tr(\theta(s_i)\,,\,s_j)\quad,\quad x_l\in E^*\,,\,s_l\in \sss\,,\,l=i,j,k.$$

Write $\rep(A,\sss)^{\aut(\sss,\tr)}$ for the  $\k$-algebra of $\aut(\sss,\tr)$-invariant
polynomial functions on the $\k$-vector space $\Hom_{\,_\k}(E, \sss)$. The symplectic
form $\om_{_\rep}=\om\otimes\tr$ makes $\rep(A,\sss)^{\aut(\sss,\tr)}$ 
into a {\it Poisson algebra}.
We have the standard Lie algebra central extension:
\begin{equation}\label{exact_R}
0\too \k \too \rep(A,\sss)^{\aut(\sss,\tr)} 
\stackrel{\delta}{\too} \derrep\bigl(\rep(A,\sss)^{\aut(\sss,\tr)}\bigr)\too 0\,,
\end{equation}
where $\dis\derrep\bigl(\rep(A,\sss)^{\aut(\sss,\tr)}\bigr)$ 
stands for the Lie algebra of derivations
of the commutative algebra $\rep(A,\sss)^{\aut(\sss,\tr)}$ respecting the Poisson bracket.
It is straightforward to check that the assignment: $\theta \mapsto
H_{_\theta}$ gives a Lie algebra splitting: 
$\dis\derrep\bigl(\rep(A,\sss)^{\aut(\sss,\tr)}\bigr)
\too\rep(A,\sss)^{\aut(\sss,\tr)}$ of the surjective morphism $\delta$ in 
the exact sequence above.

Observe next that the `infinite dimensional'
group $\aut(\A)$ acts naturally on $\Hom\palg(\A, \sss)$.
This action commutes with that of the group $\aut(\sss,\tr)$,
preserves the symplectic form $\om_{_\rep}$, but it is {\it not}
Hamiltonian, in general. That means that the induced
Lie algebra
morphism $\dis \xi: \derp\A \too 
\derrep\bigl(\rep(A,\sss)^{\aut(\sss,\tr)}\bigr)$
cannot be lifted, in general, to a Lie algebra
morphism: $\derp\A \too \rep(A,\sss)^{\aut(\sss,\tr)}\,,$
see (\ref{exact_R}).

The following result shows that
the $\aut(\A)$-action becomes  Hamiltonian after a 1-dimensional
central extension. The result below agrees also with the
philosophy advocated in [KR], saying that,
for any finite-dimensional (symmetric) $\PP$-algebra $\sss$,
`functions' on the non-commutative space
corresponding to a $\PP$-algebra $A$ should go
into genuine regular functions on the affine
algebraic variety $\Hom\palg(\A, \sss)$.

\begin{theorem}\label{Gi} There is a natural Lie algebra homomorphism
$\varphi: \R(\A) \,\to\, \rep(A,\sss)^{\aut(\sss,\tr)}$ making the following
diagram commute:
\begin{equation*}
    \xymatrix{
    {0} \ar[r] &  {\k}  \ar[r]\ar@{=}[d] &  {\enspace\R(\A)\enspace}  
\ar[r]^{\ref{central_extension}}
\ar[d]_{\varphi} &  {\enspace\derb(\A,\om)\enspace} \ar[r]\ar[d]_{\xi}&  {0}\\
{0} \ar[r]&  {\k}  \ar[r] &{\enspace\rep(A, \sss)^{\aut(\sss,\tr)}\enspace} 
\ar[r]^{\ref{exact_R}\;\;\;\quad} &
{\enspace\derrep\bigl(\rep(A,\sss)^{\aut(\sss,\tr)}\bigr)\enspace} \ar[r] &  {0}}
\end{equation*}
\end{theorem}

{\sl Proof.} Very similar to the proof of Proposition \ref{Lie}.\qed

\footnotesize{

}

\noindent
\footnotesize{
Department of Mathematics, University of Chicago, 
Chicago IL
60637, USA; $\;${\tt ginzburg@math.uchicago.edu}}

\end{document}